\documentclass[11pt,a4paper]{article}
\usepackage[margin=1in]{geometry}
\usepackage{amsmath}
\usepackage{amssymb}
\usepackage{amsthm}
\usepackage{stmaryrd}
\usepackage[utf8]{inputenc}
\usepackage[capitalize]{cleveref}
\crefname{equation}{}{}
\usepackage{tikz-cd}
\usepackage{nicefrac}
\usepackage{enumerate}
\usepackage{enumitem}
\usepackage{amsmath,amsbsy,amssymb}
\usepackage{algorithmic,algorithm}
\usepackage{multirow}
\usepackage{tikz}
\usepackage{pgfplots}
\pgfplotsset{compat=newest}
\usepackage{pgfplotstable}
\usepackage{caption}
\usepackage{subcaption}
\usepackage{overpic}
\usepackage{bm}
\usepackage[normalem]{ulem}
\usepackage[pagewise]{lineno}
\pgfplotsset{
legend image code/.code={
\draw[mark repeat=2,mark phase=2]
plot coordinates {(0cm,0cm) (0.1cm,0cm) (0.2cm,0cm)};}
}

\usepackage[normalem]{ulem}

\newtheorem{theorem}{Theorem}[section]
\newtheorem{lemma}[theorem]{Lemma}

\theoremstyle{definition}
\newtheorem{definition}[theorem]{Definition}

\theoremstyle{remark}
\newtheorem{remark}[theorem]{Remark}

\numberwithin{equation}{section}

\DeclareMathOperator*{\argmin}{arg\,min}
\newcommand{\real}{\textup{Re}}
\newcommand{\norm}[2]{\left\|#1\right\|_{#2}}

\newcommand{\normT}[1]{{\left\vert\kern-0.25ex\left\vert\kern-0.25ex\left\vert #1\right\vert\kern-0.25ex\right\vert\kern-0.25ex\right\vert}}
\newcommand{\abs}[1]{\left|#1\right|}
\newcommand{\dual}[3]{\langle #1,#2\rangle_{#3}}
\newcommand{\dualX}[2]{\dual{#1}{#2}{\XX}}

\newcommand{\XX}{\mathcal X}
\newcommand{\VV}{\mathcal V}
\newcommand{\R}{\mathbb R}
\newcommand{\C}{\mathbb C}
\newcommand{\bigO}{\mathcal O}

\usepackage{soul}
\usepackage{cancel}

\newcommand{\TT}{\mathcal{T}}
\newcommand{\MM}{\mathcal{M}}
\newcommand{\Ceff}{C_{\mathrm{eff}}}
\newcommand{\Crel}{C_{\mathrm{rel}}}
\newcommand{\Nmax}{N_{\mathrm{max}}}
\newcommand{\tol}{\mathsf{tol}}
\DeclareMathOperator*{\diam}{diam}
\DeclareMathOperator*{\refine}{refine}
\DeclareMathOperator*{\osc}{osc}

\title{Rational-approximation-based model order reduction of Helmholtz frequency response problems with adaptive finite element snapshots} 

\usepackage{authblk} 

\author[1]{Francesca Bonizzoni}
\affil[1]{MOX - Department of Mathematics, Politecnico di Milano, 20133 Milano, Italy (francesca.bonizzoni@polimi.it)}
\author[2]{Davide Pradovera}
\affil[2]{Department of Mathematics, University of Vienna, 1090 Vienna, Austria (davide.pradovera@univie.ac.at)}
\author[3]{Michele Ruggeri}
\affil[3]{Department of Mathematics and Statistics, University of Strathclyde, Glasgow G1~1XH, United Kingdom (michele.ruggeri@strath.ac.uk)}

\date{}

\begin{document}

\maketitle

\begin{abstract}
\noindent
We introduce several spatially adaptive model order reduction approaches tailored to non-coercive elliptic boundary value problems, specifically,
parametric-in-frequency Helmholtz problems. 
The offline information is computed by means of adaptive finite elements, so that each snapshot lives in a different discrete space that resolves the local singularities of the analytical solution and is adjusted to the considered frequency value.
A rational surrogate is then assembled adopting either a least-squares or an interpolatory approach, yielding  function-valued version of the standard rational interpolation method ($\mathcal{V}$-SRI) and the minimal rational interpolation method (MRI).
In the context of building an approximation for linear or quadratic functionals of the Helmholtz solution, we perform several numerical experiments to compare the proposed methodologies. Our simulations show that, for interior resonant problems (whose singularities are encoded by poles on the real axis), the spatially adaptive $\mathcal{V}$-SRI and MRI work comparably well. Instead, when dealing with exterior scattering problems, whose frequency response is mostly smooth, the $\mathcal{V}$-SRI method seems to be the best-performing one.
\end{abstract}



\section{Introduction}

Many engineering applications, e.g., in structural dynamics, geophysics, seismology, acoustics or vibro-acoustics, require the numerical approximation of solutions to time-harmonic wave propagation problems over a range of frequencies. Due to oscillations in the analytical solutions, accurate numerical approximations of frequency responses are computationally expensive and time-consuming, already for moderate frequencies. Therefore, in the multi-query context, when multiple solves of the model have to be performed for different frequency values, the ``naive'' approach of discretizing the original problem as many times as needed is usually unaffordable.

Let us denote by $u(z)$ the exact solution to the time-harmonic problem under consideration for a given wavenumber $z \in \C$. Model order reduction (MOR) methods aim at alleviating the computational cost by producing an approximation of the \emph{solution map} $z \mapsto u(z)$ or the  \emph{frequency response map} of a quantity of interest (QoI), usually given by some functional of the solution $z \mapsto y(z):=F(u(z))$.
The produced approximation (the so-called \emph{surrogate}) has to be close to the QoI and, at the same time, should be cheap to evaluate.
Customarily, MOR methods rely on a two-phase procedure. The \emph{offline} phase consists in the computation of a finite-dimensional basis of \emph{snapshots}, entailing numerical solutions of the original boundary value problem for a set of frequency values (the sample points). Based on the offline information, the surrogate for the QoI is then assembled. These steps often require a considerable computational effort. However, they are performed only once, and the surrogate is then stored for later use during the \emph{online} phase, where the surrogate is evaluated (in real-time) at any new frequency value of interest.

In the mathematical and engineering literature we can distinguish two families of MOR methods for frequency response problems. Projection-based techniques (see~\cite{BGW15}), e.g., proper orthogonal decomposition (POD), the reduced basis method (RB), and the multi-moment matching method, are widely employed and extremely powerful; however, they require access to the original problem, particularly, to the stiffness and mass matrices as well as the forcing term. In contrast, the so-called \emph{non-intrusive} MOR methods produce a surrogate relying only on the precomputed set of snapshots.
We mention, e.g., the Loewner framework~\cite{IA14}, Pad\'e-based techniques~\cite{AFR07,BNP2018}, and minimal rational interpolation (MRI) \cite{P2020, PN2020}. Non-intrusive MOR methods usually display great flexibility, since, in principle, they allow one to construct a surrogate starting from snapshots obtained via a black-box solver (e.g., commercial closed-source software). The price to pay might be a reduced accuracy (compared to intrusive methods) for a fixed set of samples.

We underline that, for frequency response problems like the Helmholtz equation (which is the main focus of this work), the analytical solution $u(z)$ can be proven to be a meromorphic function of the complex wavenumber $z$. It is therefore sensible to look for its surrogate in the class of rational $\XX$-valued maps (here $\XX=H^1(\Omega)$ with $\Omega$ being the physical domain), as all the above-mentioned MOR approaches do. Moreover, if the QoI is a linear functional of the solution field, e.g., $y(z)=F(u(z))$ with $F \in \XX'$, then it inherits a meromorphic structure from $u(z)$. The case of quadratic functionals (e.g., the energy-norm of the frequency response) is more involved, and needs to be treated separately. 

In standard MOR techniques, the snapshots of the ``truth model'' are all computed on one discretization of the considered physical domain. In the specific framework we are handling, this might represent a big drawback. Indeed, the analytical solution of the Helmholtz equation oscillates (the more so as the frequency increases), and it may exhibit local features, namely, a local resonance-type behavior, depending on the shape of the physical domain and the considered frequency values. Therefore, when a wide range of frequencies is considered, accuracy is guaranteed only by using a \emph{uniformly refined} mesh, which must be sufficiently accurate for all frequency values of interest. This can potentially entail a waste of computational resources.

In contrast, in \emph{spatially adaptive} MOR approaches, each snapshot is taken on a mesh adapted to the local features at a given parameter (i.e., wavenumber) value, and belongs to a problem-adapted finite element (FE) space that contains the ``optimal'' number of degrees of freedom (DoFs) required to obtain a certain accuracy.
Spatially adaptive MOR methods for coercive parametric partial differential equations (PDEs) have been the subject of recent mathematical research. In particular, we refer to~\cite{ASU17}, where the finite dimensional subspace for a greedy RB method consists of snapshots computed via an adaptive wavelet scheme. Adaptive FEs have been coupled with POD and greedy RB in~\cite{URL16} and \cite{Y16}, respectively.

The novelty of the present contribution resides in the introduction of spatially adaptive MOR methods for \emph{non-coercive} elliptic boundary value problems, specifically, the parametric-in-frequency Helmholtz equation. The offline information, upon which the surrogate construction relies, is a set of snapshots computed by means of an $h$-adaptive FE method (here, $h>0$ refers to the so-called \emph{mesh size}, the maximum diameter of the elements of the FE mesh).
The motivation that drives us to adaptive MOR methods is twofold. First, the use of $h$-adaptive FE snapshots, each living on a different mesh of the domain, allows to save computational resources. Second, it presents an additional advantage: while standard MOR methods aim at the approximation of the high-fidelity \emph{FE} solution, in spatially adaptive MOR methods we set as objective the approximation of the \emph{analytical} solution.

On the other hand, the adaptive approach implies intrinsic difficulties: for instance, linear combinations of snapshots cannot be easily computed. In principle, to circumvent this issue, one could express all the snapshots as elements of some common FE space. However, for adaptive mesh refinement based on newest vertex bisection (NVB) \cite{S2008,KPP13}, this would entail the construction of the so-called \emph{global} mesh overlay~\cite{S2007,CKNS08}, i.e., the smallest common refinement of all adapted meshes. In many cases, this entails a prohibitive computational effort and, more importantly, it goes against the main purpose of $h$-adaptivity. Therefore, in all the algorithms that we propose in this work, we strive to never construct the global mesh overlay. On the contrary, we will only require the evaluation of scalar products of pairs of snapshots, which is equivalent to building overlays of pairs of meshes. 

When the QoI is a scalar, one possible way to avoid dealing with snapshots on different meshes is to construct a rational surrogate for the output itself by means of standard rational interpolation (SRI) methods. Off-the-shelf rational approximation methods like AAA \cite{AAA}, the Loewner framework \cite{A05,IA14}, or vector fitting \cite{vf} can be applied to this aim. We summarize this class of approaches in \cref{sec:sri}.
    
On the other hand, the present contribution is dedicated to the development of MOR methodologies providing surrogates for function-valued QoI. With this aim:
 \begin{itemize}
       \item We extend the SRI approach to function-valued QoIs, particularly, traces and/or restrictions of the form $v(z)=u(z)|_\omega\in L^2(\omega)$ with $\omega$ being a part of the physical domain $\Omega$ or of its boundary $\partial\Omega$, or possibly even coinciding with the entire domain $\Omega$.  
    We refer to the resulting method as $\VV$-SRI. In the finite-dimensional ``vector'' case, this approach can be related to vector-valued rational approximation methods like set-valued AAA \cite{Lietaert} or the MIMO Loewner framework \cite{A05}. Still, we focus here on the infinite-dimensional case, where the quantity to approximate is function-valued.
    \item We considerably improve the MRI method \cite{P2020, PN2020} by (i) formulating it in barycentric coordinates, which allow for enhanced numerical stability properties; (ii) extending it to the $h$-adaptive FE setting. The resulting method is also referred to as $h$-adaptive MRI.
\end{itemize}

The rest of the paper is organized as follows. In \cref{sec:problems_of_interest}, the parametric-in-frequency time-harmonic wave problems of interest for our discussion (interior and scattering problems) are introduced, and their meromorphic structure is highlighted. In \cref{sec:hfem}, we recall the $h$-adaptive FE method (FEM). The core of the paper is \cref{sec:hMOR}. There, the $h$-adaptivity in the space variable is combined with the use of several rational-based MOR techniques. In \cref{sec:pod}, the $h$-adaptive POD is presented. In \cref{sec:meromorphicitysquare}, we discuss the case in which the quantity of interest is a real-quadratic functional of the analytic solution. Finally, in \cref{sec:numerical}, several numerical results are provided, to discuss the performance and highlight advantages and disadvantages of the considered approaches.

\section{Representative problems of interest}\label{sec:problems_of_interest}
\label{sec:probl_interest}

For a given complex wavenumber $k^2$, we look for $u(k^2)\in\XX:=H^1_{\Gamma_D}(\Omega)$, the weak solution of the following interior Helmholtz boundary value problem
\begin{equation}
\label{eq:bvp}
\begin{cases}
    -\Delta u(k^2) - k^2 u(k^2) = f, &\text{in }\Omega,\\
    u(k^2)=0, &\text{on }\Gamma_D\subset\partial\Omega,\\
    \partial_{\bm{\nu}} u(k^2)=g_N, &\text{on }\Gamma_N=\partial\Omega\setminus\Gamma_D.
\end{cases}
\end{equation}
In \cref{eq:bvp}, $\Omega\subset\R^d$ ($d=1,2,3$) denotes a bounded domain, whose boundary is partitioned into the (possibly empty) subsets $\Gamma_D$ and $\Gamma_N$. We denote by $\bm{\nu}$ the outward-pointing normal to $\Gamma_N$ and we assume that $f\in\XX'=H^{-1}(\Omega)$ and $g_N\in H^{-1/2}(\Gamma_N)$.

Our ultimate target is the approximation of a quantity of interest (QoI) that depends on the weak solution $u(k^2)$, for all values of the (complex) parameter $k^2$ sweeping a given range of interest $Z$. In particular, let $F:\XX\to\C$ represent a goal functional over $\XX$. Then, the target quantity for \emph{fixed} $k^2\in Z$ is the scalar
\begin{equation}\label{eq:functional}
y(k^2)=F(u(k^2))\in\C.
\end{equation}
We will consider linear as well as real-quadratic functionals. For instance, we may take
\begin{equation}\label{eq:bfunctional}
F(v)=\int_\omega v,
\end{equation}
or 
\begin{equation}\label{eq:bfunctionalq}
F(v)=\int_\omega|v|^2,
\end{equation}
(or weighted versions of them) with $\omega\subset\overline{\Omega}$ either a subdomain or a curve. Linear and real-quadratic functionals are quite often of interest in applications, representing, e.g., average displacements and vibrational energy, respectively.

For simplicity, we restrict our presentations to a $k$-independent forcing term $f$, $k$-independent Neumann datum $g_N$, and homogeneous Dirichlet boundary condition on $\Gamma_D$. 
If an inhomogeneous Dirichlet datum is given, by its lifting one falls back to a Helmholtz equation with homogeneous Dirichlet boundary conditions but with $k$-dependent right-hand side. As we showcase in our numerical examples, our discussion extends in a straightforward way to $k$-dependent forcing terms and Neumann data. However, one should note that, as the complexity of the forcing term and/or Neumann datum (particularly, with respect to $k$) increases, a surrogate of higher order -- thus entailing a higher computational cost -- is usually necessary to attain a prescribed approximation accuracy.

\begin{figure}
    \centering
     \begin{overpic}[width=0.4\textwidth]{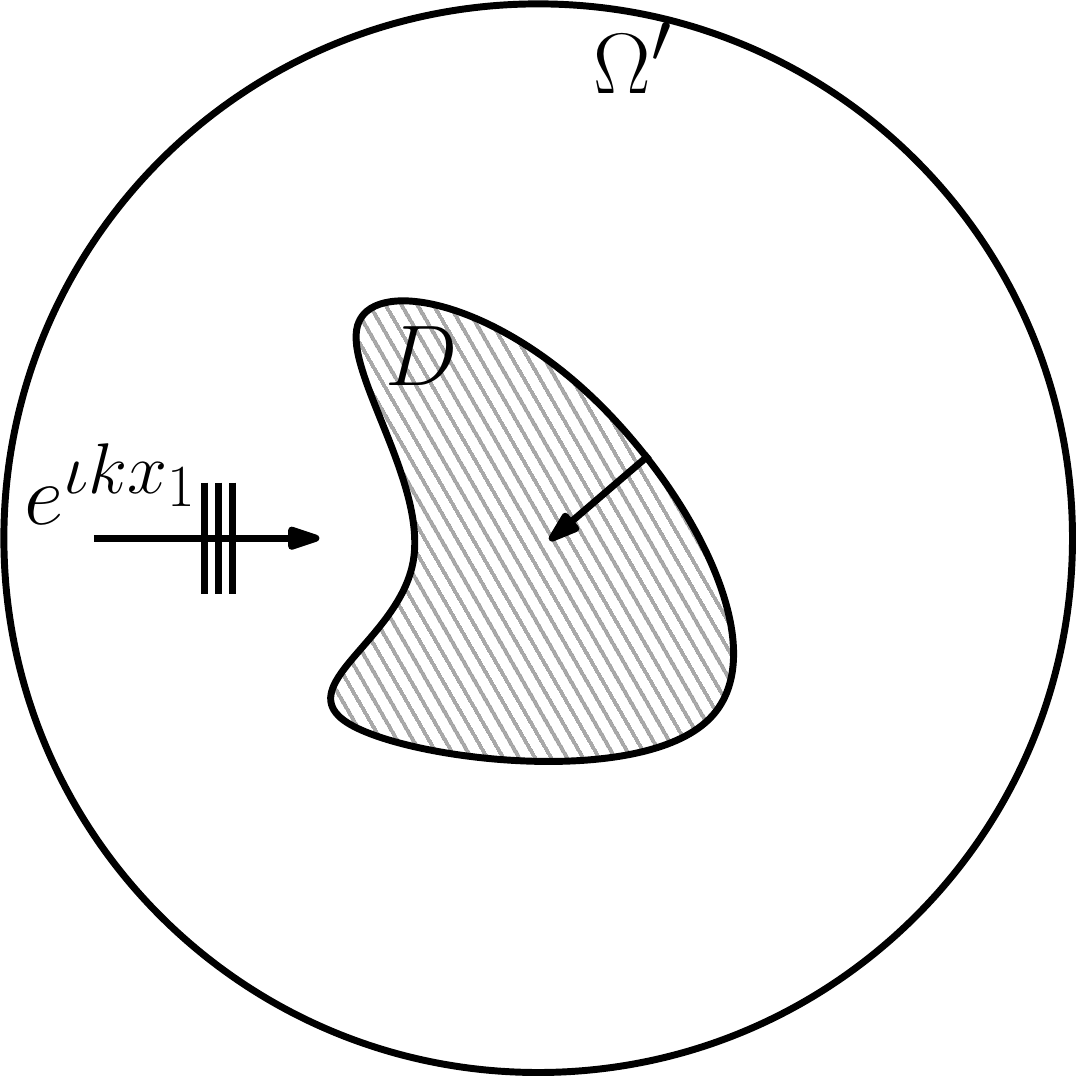}
     \put(51,55){\small${\bm \nu}$}
     \end{overpic}
    \caption{Schematic of the representative scattering problem.}
    \label{fig:domain}
\end{figure}

The interior Helmholtz problem \cref{eq:bvp} can be extended to model exterior scattering problems. Indeed, let $D\subset\R^d$ ($d=1,2,3$) denote a compact scatterer, whose surface has impedance $\zeta\in\C\cup\{\infty\}$ and inward-pointing normal $\bm{\nu}$, and define $\Omega=\R^d\setminus D$. For a given wavenumber $k\in\C$, the solution $u(k)\in\XX=H^1(\Omega)$ of
\begin{equation}\label{eq:scattering}
\begin{cases}
    -\Delta u(k) - k^2 u(k) = 0, &\text{ in } \Omega,\\
    \left(\partial_{\bm{\nu}}+\iota\frac{k}{\zeta}\right)u(k)=\left(\partial_{\bm{\nu}}+\iota\frac{k}{\zeta}\right)e^{\iota kx_1}, &\text{ on } \partial D,\\
    \lim_{|\mathbf{x}|\to\infty}|\mathbf{x}|^{(d-1)/2}\left(\partial_{|\mathbf{x}|}-\iota k\right)u(k)=0, &
\end{cases}
\end{equation}
models (in frequency domain) the wave scattered by $D$ when it is hit by the horizontal plane wave $e^{\iota kx_1}$. In practice, one can truncate $\Omega$: for instance, given $R>0$ large enough, one can set $\Omega'=\Omega\cap B(\mathbf{0},R)$, with $B(\mathbf{0},R)$ denoting a ball of radius $R$ (see \cref{fig:domain}), or $\Omega'=\Omega\cap[-R,R]^d$ (see \cref{sec:scattering}). Then, for closure, the Sommerfeld radiation condition at infinity
(the third equation in \cref{eq:scattering}) is replaced by some approximation on the boundary $\partial\Omega'\setminus\partial D$, e.g., the first-order absorbing boundary condition $\partial_{\bm{\nu}}u(k)=\iota ku(k)$.

We note that, in the exterior case, we are writing $u(k)$, as opposed to the notation $u(k^2)$ that we used for the interior case, because the wavenumber $k$ appears also with the first power in \cref{eq:scattering} (and also in the first-order absorbing boundary condition). This will allow us to use a common framework (with the argument $z$ of the solution map $u(z)$ denoting either $k^2$ or $k$) in describing the features of the solution map in the two cases (see \cref{sec:meromorphicity}).

To conclude this section, we underline that, while it is usually quite difficult to extend theoretical results from interior to exterior problems, our methods can be applied to scattering problems without the need for any modifications.

\subsection{Meromorphicity of the solution map}
\label{sec:meromorphicity}

Given the parametric-in-frequency interior boundary value problem \cref{eq:bvp}, it is natural to introduce the \emph{solution map} (or \emph{frequency response map}) $u:\C\to\XX$ that associates to each complex wavenumber squared $z=k^2$ the weak solution $u(z)$ of the corresponding Helmholtz boundary value problem.
In \cite{BNP2018} it was proven that the solution map is meromorphic over $\C$, i.e., $u$ is holomorphic over the whole complex plane, except at a set of isolated points $\Lambda=\{\lambda_i\}_{i=1}^\infty$, where it displays pole-type singularities. In particular, $\Lambda$ is the set of eigenvalues of the (minus) Laplace operator on $\Omega$ with the considered boundary conditions. Moreover, each pole has order 1, regardless of its multiplicity as Laplace eigenvalue. More precisely, we have the following expansion:
\begin{equation}\label{eq:umeromorphic}
    u(z)=\sum_{i=1}^\infty\frac{f_i\varphi_i}{\lambda_i-z},
\end{equation}
with $\varphi_i\in\XX$, $-\Delta\varphi_i=\lambda_i\varphi_i$, for all $i$, and $\{f_i\}\subset\C$ being coefficients depending on the forcing term and boundary conditions of the problem.

This result (except for the pole order) was extended to scattering problems \cref{eq:scattering} in \cite{BNPP2020}, where, in contrast, the solution map $u(z)$ is a function of the wavenumber, i.e., $z=k$.
In this case, all the poles have strictly negative imaginary part. However, a characterization as precise as in formula~\eqref{eq:umeromorphic} is not available: the poles might have order larger than 1, and the corresponding residues will not be $\XX$-orthogonal. Still, the use of rational $\XX$-valued surrogates is justified also in this setting.

Using a notation that encompasses both cases (interior and scattering problem), the frequency response map can be expressed as
\begin{equation}
    \label{eq:frm}
    u(z)=\sum_{i=1}^{\infty}\sum_{\ell=1}^{\mu_i}\frac{\zeta_{i,\ell}}{(\lambda_i-z)^\ell}=\sum_{i=1}^{\infty}\frac{\zeta_i(z)}{(\lambda_i-z)^{\mu_i}},
\end{equation}
where, for all $i$, $\lambda_i\in\C$ is a pole of $u$ with (finite) multiplicity $\mu_i$, and $\zeta_{i,\ell}\in\XX$ are its corresponding generalized residues, for $\ell=1,\ldots,\mu_i$.

The meromorphicity of the frequency response map is inherited by QoIs $y$ that are linear functionals of $u$. The poles of $y$ are a subset of the poles of $u$. In particular, $\lambda_i$ is a pole of $y$ if at least one of its residues $\zeta_{i,\ell}$ satisfies $F(\zeta_{i,\ell})\neq 0$, i.e., if $\zeta_{i,\ell}$ is not in the kernel of $F$. Real-quadratic functionals, see,
e.g., \cref{eq:bfunctionalq}, are a bit trickier, and are discussed in \cref{sec:meromorphicitysquare}. Due to the added difficulties described there, unless otherwise specified, from here onward we will assume that the QoI $y$ is a linear functional of the solution $u$.

Before proceeding, it is important to make the following observation.
\begin{remark}
Ritz- and Petrov-Galerkin projections of the problems above preserve the meromorphicity of the solution with respect to $z$. As a consequence, when, e.g., FEs are employed to discretize the PDE over a \emph{fixed} mesh of $\Omega$, the FE discrete solution is meromorphic in $z$. However, the (countably) infinite analytic poles and eigenfunctions in \cref{eq:umeromorphic} are replaced by \emph{finite dimensional} FE counterparts.
\end{remark}

\section{Solution with $h$-adaptive FEM}
\label{sec:hfem}

In this section, we describe the $h$-adaptive FEM used to compute the snapshots that serve as inputs to the proposed MOR methods.

\subsection{Problem formulation and its FE discretization}

Let us consider the boundary value problem
\begin{equation}
\label{eq:bvp_hfem}
\begin{cases}
-\Delta u - k^2 u = f&\text{in }\Omega,\\
u=0&\text{on }\Gamma_D,\\
\partial_{\bm{\nu}} u=g_N&\text{on }\Gamma_N,\\
\partial_{\bm{\nu}} u - \iota ku = g_R&\text{on }\Gamma_R,
\end{cases}
\end{equation}
where $f \in L^2(\Omega)$,
$g_N \in L^2(\Gamma_N)$
and $g_R \in L^2(\Gamma_R)$,
with $\Gamma_D$, $\Gamma_N$, $\Gamma_R$
being relatively open
and pairwise disjoint subsets of the boundary
$\partial\Omega$
forming a partition of it,
i.e.,
$\partial\Omega
= \overline{\Gamma_D} \cup \overline{\Gamma_N} \cup \overline{\Gamma_R}$.
This problem covers both~\cref{eq:bvp}
(corresponding to the case $\Gamma_R = \emptyset$)
and~\eqref{eq:scattering} with first-order absorbing boundary conditions.
The variational formulation of~\cref{eq:bvp_hfem} reads as follows:
find $u \in \XX = H^1_{\Gamma_D}(\Omega)$ such that
\begin{equation} \label{eq:weak}
b(u,v)
= \int_{\Omega} f \overline{v}
+ \int_{\Gamma_N} g_N \overline{v}
+ \int_{\Gamma_R} g_R \overline{v}
\quad
\text{for all }
v \in \XX,
\end{equation}
where the sesquilinear form $b:\XX\times\XX \to \C$ is given by
\begin{equation*}
b(w,v)
= \int_\Omega \nabla w \cdot \overline{\nabla v}
- k^2 \int_\Omega w \overline{v}
- \iota k \int_{\Gamma_R} w \overline{v}
\quad
\text{for all } w,v \in \XX.
\end{equation*}
The variational problem is well-posed for $k\in\R$, if $\Gamma_R \neq \emptyset$ or if $\Gamma_R = \emptyset$ and $k^2$ is not an eigenvalue of the Laplace operator with mixed Dirichlet-Neumann boundary conditions.

To approximate $u$,
we consider a standard $\XX$-conforming FEM.
Let $\TT_\bullet$ be a regular simplicial mesh of $\Omega$.
We assume that the partition of $\partial\Omega$ is resolved by $\TT_\bullet$
and denote by $\mathcal{E}_\bullet$ the set of facets (edges in 2D) of $\TT_\bullet$.
We consider the space of globally continuous and $\TT_\bullet$-piecewise affine functions (P1-FEM),
i.e., 
\begin{equation*}
\XX_\bullet
:= \{ v_h \in C(\overline{\Omega}) : v\vert_T \text{ is affine for all } T \in \TT_\bullet\} \cap \XX.
\end{equation*}
Then, the approximation $u_\bullet \in \XX_\bullet$ of $u$
is the solution of the following finite-dimensional variational problem:
find $u_\bullet \in \XX_\bullet$ such that
\begin{equation} \label{eq:fem}
b(u_\bullet,v_\bullet)=
\int_{\Omega} f \overline{v_\bullet}
+ \int_{\Gamma_N} g_N \overline{v_\bullet}
+ \int_{\Gamma_R} g_R \overline{v_\bullet}
\quad
\text{for all }
v_\bullet \in \XX_\bullet.
\end{equation}
It is well-known that~\cref{eq:fem} can be ill-posed even if~\cref{eq:weak} is well-posed.
However, if~\cref{eq:weak} is well-posed and $\TT_\bullet$ is sufficiently fine,
then~\cref{eq:fem} admits a unique solution;
see, e.g., \cite[Proposition~1]{BHP2017}.

\subsection{A~posteriori error estimation and adaptive algorithm}

Error estimation techniques and mesh adaptive methods are powerful tools to accelerate the convergence of FEM.
Starting with the pioneering works~\cite{br1978,BR78}, several error estimation strategies for FEM have been proposed; see, e.g., the monograph~\cite{ao2000}.
For works focused on \textsl{a~posteriori} error estimation for the Helmholtz equation, we refer, e.g., to~\cite{bisg1997a,bisg1997b,ds2013,sz2015,cev2021}.

In this work, to estimate the error between $u$ and $u_\bullet$,
we consider the classical
residual-based \textsl{a~posteriori} error estimator, i.e.,
\begin{equation*}
\eta_\bullet^2 = \sum_{T \in \TT_\bullet} \eta_\bullet(T)^2,
\end{equation*}
where,
for all $T \in \TT_\bullet$,
the \emph{computable} local refinement indicators are given by
\begin{multline} \label{eq:h-estim}
\eta_\bullet(T)^2
= h_T^2 \norm{f + k^2 u_\bullet}{L^2(T)}^2
+ h_T \norm{\llbracket\partial_{\bm{\nu}} u_\bullet\rrbracket}{L^2(\partial T \cap \Omega)}^2 \\
+ h_T \norm{g_N - \partial_{\bm{\nu}} u_\bullet}{L^2(\partial T \cap \Gamma_N)}^2
+ h_T \norm{g_R + \iota k u_\bullet - \partial_{\bm{\nu}} u_\bullet}{L^2(\partial T \cap \Gamma_R)}^2.
\end{multline}
In~\cref{eq:h-estim}, $h_T = \diam(T)$,
while $\llbracket\partial_{\bm{\nu}} u_\bullet\rrbracket$ denotes the jump of the normal derivative of $u_\bullet$ over the interface of two interior elements.
If~\cref{eq:weak} and~\cref{eq:fem} are \emph{both} well-posed,
standard arguments in \textsl{a~posteriori} error analysis
(see, e.g., \cite[Sections~2.2--2.3]{ao2000})
reveal that
$\eta_\bullet$ is \emph{reliable} and \emph{efficient} in the sense that
\begin{equation}\label{eq:Ceff}
\Crel^{-1}
\norm{\nabla(u - u_\bullet)}{L^2(\Omega)}
\le
\eta_\bullet
\le
\Ceff
\left(\norm{\nabla(u - u_\bullet)}{L^2(\Omega)} + \osc(f,g_N,g_R)\right).
\end{equation}
Here, $\Crel$ and $\Ceff$ are positive constants,
while
\begin{align*}
\osc(f,g_N,g_R)^2 = \sum_{T \in \TT_\bullet} h_T^2 \norm{f - f_T}{L^2(T)}^2 +& \sum_{e \in \mathcal{E}_\bullet \cap \Gamma_N} h_e \norm{g_N - g_{N,e}}{L^2(e)}^2\\
+& \sum_{e \in \mathcal{E}_\bullet \cap \Gamma_R} h_e \norm{g_R - g_{R,e}}{L^2(e)}^2
\end{align*}
collects the so-called data oscillations,
where $f_T$ denotes the integral mean of $f$ in $T$,
while $g_{N,e}$ (resp., $g_{R,e}$) denotes the integral mean of $g_N$ (resp., $g_R$)
in $e \in \mathcal{E}$.
The positive constants $\Crel$ and $\Ceff$ in~\cref{eq:Ceff} depend on the shape-regularity of the mesh
and on the problem data. They notably also depend on $k$.

Efficient error estimation is the fundamental ingredient to steer adaptive algorithms of the form
\begin{equation}\label{eq:hloop}
\text{SOLVE}\quad\to\quad\text{ESTIMATE}\quad\to\quad\text{MARK}\quad\to\quad\text{REFINE}.
\end{equation}
In such $h$-adaptive algorithms,
the regions of the domain characterized by large values of the local error indicators are selected for further refinement.
This mechanism automatically produces meshes that resolve the features of the solution and reduces the number of DoFs necessary to achieve a certain accuracy.
For the Helmholtz equation, this results in meshes that resolve the local singularities of the analytical solution and are adjusted to the considered frequency value.

During the last twenty years, the (optimal) convergence of adaptive FEM have been the subject of intense research; see, e.g., the recent review~\cite{cfpp2014}.
Such results guarantee that $h$-adaptive algorithms generate approximations with the desired accuracy in a finite number of iterations and that, asymptotically, the convergence has the best possible rate for the considered discretization.
We refer to the recent work~\cite{BHP2017} for the proof of optimal convergence of an adaptive FEM for compactly perturbed elliptic problems such as the Helmholtz equation.

In the present work,
for mesh refinement, we employ NVB~\cite{S2008,KPP13}.
Given a mesh $\TT_\bullet$ and a subset $\MM_\bullet \subseteq \TT_\bullet$ of marked elements,
we denote by
$\TT_\circ := \refine(\TT_\bullet,\MM_\bullet)$
the coarsest NVB refinement of $\TT_\bullet$ such that $\MM_\bullet \subseteq \TT_\bullet \setminus \TT_\circ$.
Furthermore, we assume that any mesh used to discretize $\Omega$ can be obtained
by applying a finite number of NVB refinements to a given initial mesh $\TT_0$.

Given the parameters $0 < \theta \le 1$, $\tol_h>0$, and $\Nmax > 1$, we consider the following algorithm
(see \cite[Algorithm~7]{BHP2017}):

\begin{algorithm}[H]
    \caption{$h$-adaptive FEM for the Helmholtz equation}
    \label{alg:hfem}
    \begin{algorithmic}
	\REQUIRE $z$, $\ell=0$, initial mesh $\TT_0$
	\LOOP
	\IF{\cref{eq:fem} does not admit a unique solution in $\XX_\ell$}
	    \STATE set $u_\ell := 0$, $\eta_\ell :=1$, and $\MM_\ell:=\TT_\ell$
	\ELSE
    	\STATE compute the unique solution $u_\ell \in \XX_\ell$ of~\cref{eq:fem}
	    \STATE compute the refinement indicator $\eta_\ell(T)$ in~\cref{eq:h-estim} for all $T \in \TT_\ell$
	    \IF{$\eta_\ell \le \tol_h$ \OR $\dim\XX_\ell>\Nmax$}
		    \STATE \textbf{break}
	    \ENDIF
	    \STATE determine $\MM_\ell \subseteq \TT_\ell$ of minimal cardinality satisfying\hfill($\star$)
		\begin{equation}
		\label{eq:doerfler}
		\theta \, \eta_\ell^2
		\le
		\sum_{T \in \MM_\ell} \eta_\ell(T)^2
		\end{equation}
	\ENDIF
	\STATE define $\TT_{\ell+1} := \refine(\TT_\ell, \MM_\ell)$, increase $\ell$ by 1
	\ENDLOOP
	\RETURN approximation $u_h := u_\ell$ of $u$.
    \end{algorithmic}
\end{algorithm}

Apart from the outer if-statement,
introduced in~\cite{BHP2017} to ensure that after
a finite number of iterations
the mesh becomes sufficiently fine to have well-posedness of~\cref{eq:fem},
\cref{alg:hfem} is the standard adaptive algorithm of the form \cref{eq:hloop}.
In step~($\star$),
the selection of the elements to be marked for refinement,
modulated by the parameter $0 < \theta \le 1$,
is
based on the D\"orfler marking criterion~\cref{eq:doerfler} from~\cite{D1996}.
\Cref{alg:hfem} delivers an approximation $u_h = u_h(z) \approx u(z)$ such that
either the corresponding error estimate $\eta_h$ satisfies $\eta_h \le \tol_h$
or the dimension of the underlying FE space satisfies $\dim\XX_h > \Nmax$.

\subsection{Computing bi-/sesquilinear forms across different meshes}
\label{sec:gramianoverlayed}

The construction of the surrogate via the MOR methods discussed below will require the evaluation of scalar products of pairs of snapshots $u_h(z)$ and $u_h(z')$ associated with different frequencies $z$ and $z'$. This leads to some difficulties, since, in general, they belong to the \emph{different} FE spaces $\XX_h(z)$ and $\XX_h(z')$. As such, it is necessary to assemble stiffness and mass matrices ``across'' different meshes $\TT_h(z)$ and $\TT_h(z')$.
Exploiting the fact that $\TT_h(z)$ and $\TT_h(z')$ are NVB refinements of the same fixed initial mesh $\TT_0$ and the binary tree structure of NVB, it is possible to compute the matrix entries efficiently and exactly, although the meshes are, in general, different and not even nested. For more details, we refer to \cite[Section~6.2]{BPR21}; see also \cref{sec:numerical}.

\section{Combining MOR with $h$-adaptivity}
\label{sec:hMOR}

In the present section, we describe our proposed rational-based MOR methods for the approximation of $h$-adaptive quantities. First, we recall the SRI method in barycentric coordinates, which is suited for rational approximation of scalar-valued QoIs. Then, we extend the method to function-valued QoIs. Finally, we focus on the MRI method, and we specialize it to the adaptive framework of the present work.

\subsection{Standard rational interpolation (SRI) for scalar QoIs}
\label{sec:sri}

When constructing the surrogate for a fixed QoI, a simple way to avoid dealing with snapshots on different meshes altogether is to approximate the (scalar) output $y_h(z)=F(u_h(z))$ directly. As observed in \cref{sec:meromorphicity}, linear functionals of $u$ are meromorphic, so that a rational approximation of $y_h$ is natural (see \cref{sec:thMRI} for a further discussion on this):
\begin{equation}\label{eq:outputrational}
    y_h(z)\approx y_{h,[N]}(z)=\frac{P_{[N]}(z)}{Q_{[N]}(z)},
\end{equation}
with $P_{[N]},Q_{[N]}\in\mathbb{P}_N(\C)$. 
Rational interpolation of scalar functions has been object of research for quite some time. A popular strategy entails the search of numerator and denominator of the rational approximant  by minimization of the (weighted) \emph{linearized} interpolation error at the sample points $\{z_j\}_{j=1}^S\subset\C$, i.e.,
\begin{equation}\label{eq:interpolationerror}
    \mathbb{P}_{N}(\C)\times\mathbb{P}_{N}(\C)\ni(P,Q)\;\mapsto\;\sum_{j=1}^Sw_j\left|Q(z_j)y_h(z_j)-P(z_j)\right|^2\in\R^+,
\end{equation}
with suitable weights $\{w_j\}_{j=1}^S\subset\R^+$. A normalization constraint (usually on $Q_{[N]}$) must be imposed to avoid the trivial solution $P_{[N]}=Q_{[N]}=0$. A key property of rational interpolation (in this form) is that at least $2N+1$ samples are necessary. 
We note that, in the most general formulation of rational approximation, the numerator and denominator might have different polynomial degree. We ignore this here by allowing them to have defective degrees.

In the interest of numerical stability, a very useful representation of the polynomials $P_{[N]}$ and $Q_{[N]}$ can obtained in the so-called barycentric form:
\begin{equation}\label{eq:interpolationbarycentric}
    P_{[N]}(z)=\prod_{i=0}^N(z-\zeta_i)\sum_{j=0}^N\frac{p_j}{z-\zeta_j}\quad\text{and}\quad Q_{[N]}(z)=\prod_{i=0}^N(z-\zeta_i)\sum_{j=0}^N\frac{q_j}{z-\zeta_j},
\end{equation}
where $\{\zeta_j\}_{j=0}^N\subset\C$ are given distinct \emph{support points}. Henceforth, we will refer to $\pi(z)=\prod_{i=0}^N(z-\zeta_i)$ as \emph{nodal polynomial}. Whenever a support point $\zeta_i$ coincides with a sample point $z_j$, interpolation at that point is guaranteed by setting $p_i=q_iy_h(\zeta_i)$, and the $j$-th addend can (and should) be removed from \cref{eq:interpolationerror}. With this choice of basis, it is common to weigh the interpolation problem by nodal polynomial values: $w_j=|\pi(z_j)|^{-2}$. This allows to express \cref{eq:interpolationerror} in an extremely simple form.

\begin{definition}[SRI]
\label{def:sri}
    Let (distinct) sample points $\{z_j\}_{j=1}^S$ and (distinct) support points $\{\zeta_i\}_{i=0}^N$ be given, with $S\geq 2N+1$. The (barycentric) standard rational interpolant (SRI) of type $[N]$ of $y_h$ based on such sample and support points is a ratio $y_{h,[N]}^\textup{SRI}=P_{[N]}^\textup{SRI}\big/Q_{[N]}^\textup{SRI}$, with $P_{[N]}^\textup{SRI}$ and $Q_{[N]}^\textup{SRI}$ of the form~\cref{eq:interpolationbarycentric}. The coefficients $\{p_i\}_{i=0}^N\cup\{q_i\}_{i=0}^N\subset\C$ are chosen so that they minimize
    \begin{equation}\label{eq:interpolationerrorbary}
        \sum_{\substack{j=1\\z_j\notin\{\zeta_i\}_{i=0}^N}}^S\left|\sum_{i=0}^N\frac{q_iy_h(z_j)-p_i}{z_j-\zeta_i}\right|^2
    \end{equation}
    under the constraints
    \begin{equation}\label{eq:interpolationconstraints}
    \begin{cases}
        \sum_{i=0}^N\left|q_i\right|^2=1,&\\
        p_i=q_iy_h(\zeta_i) & \text{for all }i=0,\ldots,N\text{ such that }\zeta_i\in\{z_j\}_{j=1}^S.
    \end{cases}
    \end{equation}
    (For more information on the barycentric rational form, see, e.g., \cite{AAA,Klein}, and the references therein.)
\end{definition}

Similar (sometimes, equivalent) definitions yield some of the state-of-the-art algorithms for rational approximation, e.g., AAA \cite{AAA}, the Loewner framework \cite{A05}, and vector fitting \cite{vf}.

\begin{remark}
    Let $i$ be such that $q_i\neq 0$. It is not difficult to see from \cref{eq:interpolationbarycentric} that $y_{h,[N]}^\textup{SRI}(\zeta_i)=p_i/q_i$. In particular, this means that, if $\zeta_i\in\{z_j\}_{j=1}^S$ and $q_i\neq 0$, then we have interpolation of the target function at $\zeta_i$: $y_{h,[N]}^\textup{SRI}(\zeta_i)=y_h(\zeta_i)$. This provides an intuitive justification for choosing the support points as a subset of the sample points $\{\zeta_i\}_{i=0}^N\subset\{z_j\}_{j=1}^S$, as done, e.g., in the AAA method \cite{AAA}.
\end{remark}

\begin{remark}
    In \cref{def:sri}, we have imposed a normalization on the Euclidean norm of the coefficients of the denominator $Q_{[N]}^\textup{SRI}$ to exclude the trivial solution $P_{[N]}^\textup{SRI}=Q_{[N]}^\textup{SRI}=0$. This choice has been observed to be more numerically robust than the standard one, based on the normalization of a single coefficient of the denominator, and is the \emph{de facto} standard in rational approximation; see, e.g., \cite{Gonnet2011,AAA}.
\end{remark}

Before proceeding, we outline here a practical algorithm for building SRIs. In the most common \emph{collocation} case ($\{\zeta_i\}_{i=0}^N\subset\{z_j\}_{j=1}^S$), the coefficients $q_i$ can be found from the SVD of the Loewner matrix $G\in\C^{(S-N-1)\times(N+1)}$, defined entry-wise as
\begin{equation}\label{eq:sriloewner}
G_{(j-N-2)i}=\frac{y_h(z_j)-y_h(z_{i+1})}{z_j-z_{i+1}}\quad\text{for }j=N+2,\ldots,S\text{ and }i=0,\ldots,N.
\end{equation}
For more details, we refer either to~\cite{AAA,Gosea21,IA14}, to \cref{alg:vsri} in the next section, or to \cref{rem:latsrisri}.

\subsection{SRI for function-valued quantities}
\label{sec:gsri}

The construction above can be easily generalized to vector quantities of interest, e.g., $\mathbf{y}_h(z)=\left(y^1_h(z),\ldots,y^O_h(z)\right)^\top$. To this aim, two changes are necessary: first, while the denominator remains $\C$-valued, the numerator must become $\C^O$-valued, i.e., $\{\mathbf{p}_i\}_{i=0}^N\subset\C^O$; second, the absolute value in \cref{eq:interpolationerror,eq:interpolationerrorbary} must be replaced by the Euclidean norm over $\C^O$, yielding
\begin{equation}\label{eq:interpolationoptimalityvector}
    \sum_{\substack{j=1\\z_j\notin\{\zeta_i\}_{i=0}^N}}^S\norm{\sum_{i=0}^N\frac{q_i\mathbf{y}_h(z_j)-\mathbf{p}_i}{z_j-\zeta_i}}{\C^O}^2=\sum_{\substack{j=1\\z_j\notin\{\zeta_i\}_{i=0}^N}}^S\sum_{l=1}^O\left|\sum_{i=0}^N\frac{q_iy_h^l(z_j)-(\mathbf{p}_i)_l}{z_j-\zeta_i}\right|^2.
\end{equation}
We note that, in the collocation case ($\{\zeta_i\}_{i=0}^N\subset\{z_j\}_{j=1}^S$), this strategy is closely related (in fact, mostly equivalent) to the ``set-valued AAA'' and ``fast-AAA'' algorithms introduced in \cite{Lietaert} and \cite{Hochman}, respectively. An SVD-based solution is possible here as well, but we skip the details, since we provide them below in the more general framework of function-valued rational approximation, which includes this as special case. Alternatively, the MIMO Loewner framework \cite{Gosea21,IA14} solves this kind of approximation problem by tangential interpolation, effectively recasting it as modified scalar problem.

More interesting and definitely less trivial is the extension to ``infinite-dimensional'' function-valued QoI. In the following, taking inspiration from \cref{eq:bfunctional,eq:bfunctionalq}, we focus on the case
\begin{equation}\label{eq:bfunctionalinfinite}
    v(z)=u(z)|_\omega\in L^2(\omega)=:\mathcal{V}.
\end{equation}
However, as we describe in more detail in~\cref{sec:thMRI}, we note that our discussion applies also to the case $\omega=\Omega$, i.e., $v=u$.
We obtain the following definition by extending \cref{eq:interpolationoptimalityvector} from $\C^O$ to $\mathcal{V}$.
\begin{definition}[$\mathcal{V}$-SRI]\label{def:xsri}
    Let (distinct) sample points $\{z_j\}_{j=1}^S$ and (distinct) support points $\{\zeta_i\}_{i=0}^N$ be given, with $S\geq 2N+1$. The (barycentric) $\mathcal{V}$-SRI of type $[N]$ of $v_h:\C\to\mathcal{V}$ based on such sample and support points is a ratio $v_{h,[N]}^{\mathcal{V}\textup{-SRI}}=P_{[N]}^{\mathcal{V}\textup{-SRI}}\big/Q_{[N]}^{\mathcal{V}\textup{-SRI}}$, with $P_{[N]}^{\mathcal{V}\textup{-SRI}}$ and $Q_{[N]}^{\mathcal{V}\textup{-SRI}}$ being of the form~\cref{eq:interpolationbarycentric}. The coefficients $\{p_i\}_{i=0}^N\subset\mathcal{V}$ and $\{q_i\}_{i=0}^N\subset\C$ are chosen so that they minimize
    \begin{equation}\label{eq:interpolationgammatarget}
        \sum_{\substack{j=1\\z_j\notin\{\zeta_i\}_{i=0}^N}}^S\left\|\sum_{i=0}^N\frac{q_iv_h(z_j)-p_i}{z_j-\zeta_i}\right\|_{\VV}^2
    \end{equation}
    under the constraints
    \begin{equation}\label{eq:interpolationgammaconstraints}
    \begin{cases}
        \sum_{i=0}^N\left|q_i\right|^2=1,&\\
        p_i=q_iv_h(\zeta_i) & \text{for all }i=0,\ldots,N\text{ such that }\zeta_i\in\{z_j\}_{j=1}^S.
    \end{cases}
    \end{equation}
\end{definition}

\begin{remark}
\label{rem:xsristar}
    In \cref{eq:bfunctionalinfinite}, we have set $\VV=L^2(\omega)$ with the objective of approximating some trace of $u_h$. However, \cref{def:xsri} may be applied also in more general spaces $\mathcal{V}$. For instance, we may set $\mathcal{V}=\XX$ (the space where the solution $u$ lives) to define $\XX$-SRI.
\end{remark}

It is crucial to observe that the following property (whose proof is deferred to \cref{sec:traceinterpolationspan}) holds true.
\begin{lemma}\label{lem:tracepspan}
    The coefficients of the $\VV$-SRI numerator $P_{[N]}^{\VV\textup{-SRI}}$ belong to the span of the snapshots:  $\{p_i\}_{i=0}^N\subset\textup{span}\{v_h(z_j)\}_{j=1}^S$.
\end{lemma}

Accordingly, there exists a matrix $\mathring{P}\in\C^{S\times(N+1)}$ of expansion coefficients:
\begin{equation}\label{eq:tracepexpansion}
    p_i=\sum_{j=1}^S\mathring{P}_{ji}v_h(z_j)\quad\text{for }i=0,\ldots,N,
\end{equation}
and the $\VV$-SRI admits the expansion
\begin{equation}\label{eq:traceinterpolatorsurrogate}
    v_{h,[N]}^{\VV\textup{-SRI}}(z)=\left(\sum_{i=0}^N\frac{\sum_{j=1}^S\mathring{P}_{ji}v_h(z_j)}{z-\zeta_i}\right)\bigg/\left(\sum_{i=0}^N\frac{q_i}{z-\zeta_i}\right).
\end{equation}

An algorithm for computing $\VV$-SRIs can be obtained as an extension of that for SRI. Notably, as we detail in \cref{sec:traceinterpolation}, the minimization of \cref{eq:interpolationgammatarget} can be carried out through an SVD-like procedure in the ``vectorized'' $\VV\otimes\C^S$ metric. As a practical way to do this, we introduce a $\VV$-orthonormalization step: given $\{v_h(z_j)\}_{j=1}^S\subset\VV$, we compute a $\VV$-orthonormal ($\langle\psi_j,\psi_{j'}\rangle_{\VV}=\delta_{jj'}$) basis $\{\psi_{j'}\}_{j'=1}^T\subset\VV$ and a $T\times S$ matrix $R$ such that
\begin{equation}\label{eq:householder}
    v_h(z_j)=\sum_{j'=1}^TR_{j'j}\psi_{j'}\qquad\text{for }j=1,\ldots,S.
\end{equation}
(The value $T\leq S$ is the rank of the snapshots, i.e., the dimension of their span as a subspace of $\VV$.) The basis $\{\psi_{j'}\}_{j'=1}^T$ can be constructed, e.g., by Householder triangularization of the snapshot \emph{quasi}-matrix $[v_h(z_1)|\cdots|v_h(z_S)]$, see \cite{Householder}. Alternatively, one can obtain $R$ from a Cholesky decomposition of a suitable Gramian matrix, cf.\ the proof of \cref{lem:latsri}.

\begin{algorithm}
    \caption{$\VV$-SRI}
    \label{alg:vsri}
    \begin{algorithmic}
        \REQUIRE sample points $\{z_j\}_{j=1}^S\subset\C$ and denominator degree $N\leq\frac{S-1}{2}$
        \REQUIRE target function $v_h:\C\to\VV$
        \STATE compute $v_h(z_1),\ldots,v_h(z_S)$
		\STATE orthonormalize the snapshots to obtain $R\in\C^{T\times S}$, see \cref{eq:householder}
		\STATE assemble the (vectorized) Loewner matrix $G$ defined in \cref{eq:vsriloewner}
		\STATE compute the SVD of $G=U\Sigma V^H$, with $\Sigma\in\C^{T(S-N-1)\times(N+1)}$
		\STATE define $(q_0,\ldots,q_N)^\top=V_{:N}$, the last column of $V$
		\RETURN surrogate $\left(\sum_{i=0}^N\frac{q_iv_h(z_{i+1})}{\cdot\ -z_{i+1}}\right)/\left(\sum_{i=0}^N\frac{q_i}{\cdot\ -z_{i+1}}\right)$
    \end{algorithmic}
\end{algorithm}

We summarize the overall method in \cref{alg:vsri}. For simplicity, we only consider the framework where the support points are a subset of the sample points, with $\zeta_i=z_{i+1}$ for $i=0,\ldots,N$. We note that the algorithm relies on the Loewner matrix $G\in\C^{T(S-N-1)\times(N+1)}$, defined entry-wise as
\begin{equation}\label{eq:vsriloewner}
G_{(T(j-N-2)+j')i}=\frac{R_{j'j}-R_{j'(i+1)}}{z_j-z_{i+1}}
\end{equation}
for $j=N+2,\ldots,S$, $j'=1,\ldots,T$, and $i=0,\ldots,N$, which is just a vectorized version of \cref{eq:sriloewner}.

The general case ($\{\zeta_i\}_{i=0}^N\not\subset\{z_j\}_{j=1}^S$) allows for a similar algorithm, where, however, the coefficients of the numerator are not necessarily scalar multiples of snapshots. We provide a strategy to compute $\VV$-SRIs in \cref{sec:traceinterpolation}.

\subsection{Minimal rational interpolation (MRI) for state estimation}
\label{sec:mri}

As we will motivate in~\cref{sec:thMRI}, the approximation of the full state $u(z)$ is a topic that deserves particular attention \emph{per se}. Due to its importance, several MOR techniques start by computing a surrogate for the state $u$ and only afterwards derive from it an approximation for the quantity of interest $y$. Among these, we can find projective approaches (see \cref{sec:pod} and the references therein), least-squares Pad\'e-based approaches (see \cite{BNPP2020fast}) and the previously mentioned $\XX$-SRI. This last technique has been further developed in~\cite{P2020} under the name of \emph{minimal rational interpolation} (MRI), with the objective of reducing the number of snapshots $S$ needed to achieve a certain rational type $[N]$ (or, vice versa, maximizing the rational type for a given number of snapshots). In particular, the MRI methodology allows constructing a rational approximant of type $[S-1]$, given only $S$ samples of $u$ (as opposed to the $2S-1$ that would be necessary in, e.g., SRI). 
In the following paragraphs, we not only extend MRI to the $h$-adaptive setting but we also formulate it in barycentric coordinates.

\begin{definition}[MRI]
    The (barycentric) MRI of $u_h$ based on (distinct) sample points $\{z_j\}_{j=1}^S$ is the ratio
    \begin{equation*}
        u_{h,[S-1]}^\textup{MRI}(z)=\frac{P_{[S-1]}^\textup{MRI}(z)}{Q_{[S-1]}^\textup{MRI}(z)}=\left(\pi(z)\sum_{j=1}^S\frac{q_{j-1}u_h(z_j)}{z-z_j}\right)\Bigg/\left(\pi(z)\sum_{j=1}^S\frac{q_{j-1}}{z-z_j}\right),
    \end{equation*}
    $\pi$ being the nodal polynomial $\pi(z)=\prod_{j=1}^S(z-z_j)$, where $\{q_i\}_{i=0}^{S-1}\subset\C$ minimizes
    \begin{equation}\label{eq:mrioptimality}
        \left\|\sum_{j=1}^Sq_{j-1}u_h(z_j)\right\|_\XX^2
    \end{equation}
    under the constraint $\sum_{i=0}^{S-1}\left|q_i\right|^2=1$.
\end{definition}

\begin{remark}\label{rem:mridominant}
As seen in the previous sections, the choice of coefficients of $P_{[S-1]}^\textup{MRI}$ ensures interpolation of $u_h$ at all support points, which here coincide with \emph{all} the sample points. Moreover, we observe that the target quantity \cref{eq:mrioptimality} corresponds to the $\XX$-norm of the leading coefficient of the numerator:
\begin{equation*}
    \frac{1}{(S-1)!}\frac{\textup{d}^{S-1}P_{[S-1]}^\textup{MRI}}{\textup{d}z^{S-1}}=\sum_{j=1}^S\frac{1}{(S-1)!}\frac{\textup{d}^{S-1}}{\textup{d}z^{S-1}}\left(\frac{\pi(z)}{z-z_j}\right)q_{j-1}u_h(z_j)=\sum_{j=1}^Sq_{j-1}u_h(z_j).
\end{equation*}
\end{remark}

It is interesting to note that the MRI approximant may be cast in the form \cref{eq:traceinterpolatorsurrogate} used for $\XX$-SRI (or $L^2(\Omega)$-SRI), where the coefficient matrix $\mathring{P}\in\C^{S\times S}$ is given by:
\begin{equation*}
    \mathring{P}_{ji}=\delta_{(i+1)j}q_i\quad\textup{for }i=0,\ldots,S-1\text{ and }j=1,\ldots,S.
\end{equation*}
In fact, MRI could be interpreted as an extension of $\XX$-SRI with an ``unnatural'' choice $S=N+1$.

The characterization of the denominator coefficients $\mathbf{q}=(q_0,\ldots,q_{S-1})^\top$ in the definition of MRI can be equivalently stated as: $\mathbf{q}$ is a minimal eigenvector of the snapshot Gramian
\begin{equation}\label{eq:mrisnapgram}
    G_h^{(u)}=
    \begin{bmatrix}
    \|u_h(z_1)\|_\XX^2 & \cdots & \langle u_h(z_S),u_h(z_1)\rangle_\XX \\
    \vdots & \ddots & \vdots\\
    \langle u_h(z_1),u_h(z_S)\rangle_\XX & \cdots & \|u_h(z_S)\|_\XX^2
    \end{bmatrix}
    \in\C^{S\times S},
\end{equation}
cf.\ \cref{eq:interpolationsnapgram}. This helps in designing a numerical strategy to compute the MRI surrogate, which we report in \cref{alg:mri}.

\begin{algorithm}
    \caption{MRI}
    \label{alg:mri}
    \begin{algorithmic}
        \REQUIRE sample points $\{z_j\}_{j=1}^S\subset\C$ and target function $u_h:\C\to\XX$
        \STATE compute $u_h(z_1),\ldots,u_h(z_S)$
		\STATE assemble the snapshot Gramian $G_h^{(u)}$ defined in \cref{eq:mrisnapgram}
		\STATE compute the SVD of $G_h^{(u)}=V\Sigma V^H$
		\STATE define $(q_0,\ldots,q_{S-1})^\top=V_{:(S-1)}$, the last column of $V$
		\RETURN surrogate $\left(\sum_{j=1}^S\frac{q_{j-1}u_h(z_j)}{\cdot\ -z_j}\right)/\left(\sum_{j=1}^S\frac{q_{j-1}}{\cdot\ -z_j}\right)$
    \end{algorithmic}
\end{algorithm}

\subsection{Some considerations on $\VV$-SRI and MRI in the $h$-adaptive case}\label{sec:thMRI}

In the $h$-adaptive setting, the samples of the quantities $v_h$ (see~\cref{sec:gsri}) and $u_h$ live on potentially different meshes. This represents an additional difficulty with respect to standard MOR (where the snapshots are all computed using the same mesh), which we must take it into account when defining, building, and evaluating the surrogate. For instance, evaluating the $L^2(\omega)$-SRI surrogate requires combining snapshots over all the (potentially different) meshes discretizing $\omega$, resulting in an overly high online cost due to a global mesh overlay. Unfortunately, this issue is intrinsic to the task of approximating a non-local quantity and cannot be easily mitigated. On the other hand, once $v_{h,[N]}^{L^2(\omega)\textup{-SRI}}\approx u|_\omega$ has been computed, it is easy to derive \textsl{a~posteriori} (without training a new surrogate) surrogates for arbitrary linear functionals of $u(z)$, provided their support lies in $\omega$. One possible example is
\begin{equation*}
    y(z)=F(u(z))=\int_{\omega}wu(z),
\end{equation*}
with $w:\omega\to\C$ some weight function. In such cases, it suffices to extract the scalar samples $\{y_h(z_j)\}_{j=1}^S$ from the infinite-dimensional ones $\{v_h(z_j)\}_{j=1}^S$, and then replace $v_h$ by $y_h$ in \cref{eq:traceinterpolatorsurrogate}. Since the support of $F$ is contained in $\omega$, this can be done without using the full state $u_h(z)$.

In view of this, computing the $\XX$-SRI (or the MRI, or even the $L^2(\Omega)$-SRI) of $u(z)$ is quite powerful, since it allows to obtain \textsl{a~posteriori} a surrogate for any linear functional of $u(z)$, without the need to repeat the training phase. As we will see in our numerical experiments in \cref{sec:numerical}, the main price to pay for this flexibility is the offline computational cost needed to build the snapshot Gramian matrix \cref{eq:interpolationsnapgram}, which is necessary to find $R$. Indeed, each of its $\frac{N(N-1)}2$ independent non-diagonal entries requires building a (potentially different) pairwise mesh overlay over the whole $\Omega$, cf.\ \cref{sec:gramianoverlayed}.

Before proceeding further, we wish to mention that, while MRI has the advantage of using the samples optimally (and of being extremely simple to implement), it has one main drawback: its theoretical foundations hold only when the target of the surrogate admits a simple partial fraction decomposition whose residues are linearly independent. More explicitly, the following must hold true for the theory in \cite{P2020} (particularly, the convergence of $u_{[S-1]}^\textup{MRI}$ to $u$ in the $\XX$-norm as $S\to\infty$) to apply:
\begin{enumerate}
    \item $u$ is of the form \cref{eq:frm}, with simple poles, i.e., $\mu_i=1$ for all $i$; note that the expansion has to hold only for $z$ in a (large enough) neighborhood of the target wavenumber range $Z$;
    \item the $S$ sample points $\{z_j\}_{j=1}^S$ are judiciously placed on the wavenumber range of interest (good Lebesgue constant);
    \item the $N+1=S$ most relevant residues (relevance is determined according to the Green's potential of $Z$, see \cite{P2020}) should be linearly independent.
\end{enumerate}
The first point is quite relevant to our discussion, since, in our setting, $u$ is of the required form, see \cref{eq:umeromorphic}, but $u_h$ is not. A way to formalize this issue is the following: we take as approximation target the analytic solution $u$, whose snapshots are, however, affected by the noise $u_h-u$. From this viewpoint, an analysis of the properties of $h$-MRI could be based on the stability of MRI approximation with respect to noise. A theoretical discussion on this issue is planned to appear in a forthcoming paper. In the present work, we only observe the (sometimes adverse) effects of the $h$-FEM noise in our numerical experiments in \cref{sec:numerical}. Of course, due to the adaptive nature of the snapshots, these noise-related issues are intrinsic to the problem at hand. As such, they may affect all the presented methods. Still, the LS-based formulations of SRI and $\VV$-SRI have a regularizing effect that is missing in MRI.

\section{Projective MOR alternatives}
\label{sec:pod}

In the literature, many MOR approaches have been proposed to approximate (functionals of) the solution of the Helmholtz equation. So-called \emph{projective} methods build an approximation of the solution $u$, rather than of $y$ directly, through a restriction of the state problem (the Helmholtz equation) onto a suitably chosen subspace of $\XX$, called ``reduced space''. Then the surrogate for $y$ is extracted from that of $u$, usually at small computational cost. The main feature characterizing a projective MOR technique is how such subspace is selected.

In POD, $u$ is evaluated at some prescribed parametric values, building a collection of snapshots. A principal components analysis of the snapshots is performed (e.g., by a generalized SVD), allowing the identification of their dominant modes. The reduced space is then defined as the span of such modes. We refer to~\cite{QMN15} for an extensive discussion of POD.
Then, the restriction of the original problem onto the reduced space is performed by (Petrov-)Galerkin projection onto $\widetilde{\XX}=\textup{span}\{\psi_j\}_j$, with $\{\psi_j\}_j$ the \emph{reduced basis} of choice. Notably, this relies on a so-called \emph{affine} decomposition of the original problem: Find $u(z)\in\XX$ such that
\begin{equation}\label{eq:affine}
    \sum_{\ell=1}^{n_A}\phi_\ell(z)A_\ell u(z)=\sum_{\ell=1}^{n_b}\psi_\ell(z)b_\ell,\quad\phi_\ell,\psi_\ell:\C\to\C,\ A_\ell:\XX\to\XX,\ b_\ell\in\XX,
\end{equation}
gets projected onto $\widetilde{\XX}$ as: Find $\widetilde{u}(z)=\sum_ju_j(z)\psi_j\in\widetilde{\XX}$ such that
\begin{equation}\label{eq:affinereduced}
    \sum_j\sum_{\ell=1}^{n_A}\phi_\ell(z)\dualX{A_\ell\psi_j}{\psi_i} u_j(z)=\sum_{\ell=1}^{n_b}\psi_\ell(z)\dualX{b_\ell}{\psi_i}\;\quad\forall i.
\end{equation}
Obviously, this procedure requires access to each term $A_\ell$ and $b_\ell$ of the decomposition\footnote{As is evident from \cref{eq:affinereduced}, access to the (discretized) operator $A_\ell$ is usually necessary only through ``matrix-vector multiplies'', i.e., it must be possible to query $A_\ell\psi$ for an arbitrary $\psi\in\XX$.}. Strategies are available to compute affine approximations of non-affine problems (see, e.g., \cite{QR14}). In our $h$-adaptive setting, additional difficulties arise, since snapshots on different meshes are incompatible at the discrete level. Consequently, \emph{each} of the operators $A_\ell$ appearing in \cref{eq:affine} must undergo a treatment similar to the one described in \cref{sec:gramianoverlayed}, see, e.g.,~\cite{URL16}.

In the scope of this paper, we can highlight three main differences between MRI and POD:
\begin{itemize}
    \item Approximation quality: for a fixed number of snapshots, Galerkin projection identifies a quasi-optimal approximation on their span, whereas for general problems no guarantee on the accuracy of MRI is available. In specific cases, e.g., for all Helmholtz problems of the form \cref{eq:bvp}, a theory for MRI can be derived, leading to quasi-optimality guarantees \cite{P2020}. In such situations, the two approaches are often very comparable, cf.\ our numerical results in \cref{sec:numerical}. In comparison, SRI or $\VV$-SRI usually require twice as many snapshots to achieve a certain approximation order\footnote{\label{ftn:order}By ``approximation order'', for rational interpolation we mean the degree of the surrogate denominator plus one. For POD/RB, we mean the size of the reduced basis.}.
    \item Unfavorable problems: being based only on snapshots of $u$, MRI is somewhat oblivious to the structure of the underlying problem, which has $u$ as solution. As a consequence, some problems are inherently \emph{not} amenable to approximation by MRI, while POD, being extremely structure-aware, manages to work well for them. This is an intrinsic difficulty for non-intrusive methods.
    \item Affinity requirements: the necessity for an affine (or affinely approximable by, e.g., EIM \cite{QMN15}) problem structure limits the applicability of POD. This is particularly relevant for scattering problems like \cref{eq:scattering}. Indeed, an affine decomposition of the impinging wave $e^{\iota kx_1}$ must necessarily include many terms, especially if the wavenumber range is large, thus hindering either online efficiency or surrogate accuracy (or both).
\end{itemize}

A schematic comparison of the presented algorithms, and of their computational costs, is included in \cref{tab:compare}.

\begin{table}[p]
    \begin{center}
    \begin{small}
    \renewcommand{\arraystretch}{1.5}
    \begin{tabular}{c|c|c|c|c|}
     & SRI & $\VV$-SRI & MRI & POD\\
    \hline
    \hline
    \multirow{3}{*}{\rotatebox{90}{Target}} & \multicolumn{4}{c|}{Linear functional} \\
    \cline{2-5}
     & & \multicolumn{3}{c|}{Real-quadratic functional} \\
    \cline{2-5}
     & & \multicolumn{3}{c|}{Restriction to sub-domain or (sub-)boundary} \\
    \hline
    \hline
    \multirow{2}{*}{\rotatebox{90}{Problem}} & \multicolumn{4}{c|}{Affine and meromorphic}\\
    \cline{2-5}
     & \multicolumn{3}{c|}{Non-affine} & Non-meromorphic\\
    \hline
    \hline
    \multirow{12}{*}{\rotatebox{90}{Offline phase}} & \multicolumn{4}{c|}{Compute snapshots $u_h$ and meshes, extract outputs $y_h$} \\
     & \multicolumn{4}{c|}{$\bigO(Sa_h)$} \\
    \cline{2-5}
     & Build $G_h^{(y)}$ & Build $G_h^{(v)}$ & Build $G_h^{(u)}$ & Project affine \\
     & as in \cref{eq:interpolationsnapgram} & as in \cref{eq:interpolationsnapgram} & as in \cref{eq:mrisnapgram} & LHS terms \\
     & $\bigO(S^2)$ & $\bigO(S^2n_{h,\omega}^2)$ & $\bigO(S^2n_h^2)$ & $\bigO(S^2n_An_h^2)$ \\
    \cline{2-5}
     & Build Loewner & Build vectorized & & Project affine\\
     & matrix & Loewner matrix & & RHS terms \\
     & $\bigO(S^2)$ & $\bigO(S^3)$ & & $\bigO(Sn_bn_h)$ \\
    \cline{2-5}
     & \multicolumn{3}{c|}{Build rational approximant by SVD} & \\
     & \multicolumn{3}{c|}{$\bigO(S^3)$} & \\
    \cline{2-5}
     & & \multicolumn{2}{c|}{Extract approximation of $y$} & \\
     & & \multicolumn{2}{c|}{$\bigO(S^2)$} & \\
    \hline 
    \hline 
    \multirow{6}{*}{\rotatebox{90}{Online phase}} & \multicolumn{3}{c|}{ } & Build surrogate \\
     & \multicolumn{3}{c|}{ } & $\bigO(S^2n_A+Sn_b)$ \\
    \cline{5-5}
     & \multicolumn{3}{c|}{Evaluate rational approximant} & Solve surrogate \\
     & \multicolumn{3}{c|}{$\bigO(S)$} & $\bigO(S^3)$ \\
    \cline{5-5}
     & \multicolumn{3}{c|}{ } & Apply functional \\
     & \multicolumn{3}{c|}{ } & $\bigO(S)$ \\
    \hline
    \end{tabular}
    \renewcommand{\arraystretch}{1}
    \begin{tabular}{llr}
    \multicolumn{3}{l}{Legend:} \\
    & $S=\,$ number of snapshots taken & \\
    & $n_h=\,$ representative number of FE DoFs of $u_h$ & \\
    & $n_{h,\omega}=\,$ representative number of FE DoFs of $u_h|_\omega$ & ($\sim n_h^{1-1/d}$)\\
    & $a_h=\,$ representative complexity for computing $u_h$ & ($\lesssim n_h^2$)\\
    & $n_A,n_b=\,$ number of affine terms in \cref{eq:affine} & \\
    \end{tabular}
    \end{small}
    \end{center}
    \vspace{.25cm}
    \caption{\label{tab:compare}Summary of described MOR algorithms (with snapshot locations fixed \emph{a priori}) and corresponding complexities. For simplicity, we assume $N\sim S/2$ for SRI and $\VV$-SRI, leading to ``optimal'' snapshot usage in such methods. In the offline-online steps, for simplicity we assume that $y$ is a linear functional with support in $\omega\subset\partial\Omega$. In some cases, the online phase for POD may be reduced to $\bigO(Sn_b)$ complexity at $\bigO(S^3)$ offline cost.}
\end{table}

\section{Rational interpolation for real-quadratic outputs}\label{sec:meromorphicitysquare}

The case of real-quadratic functionals differs from that of linear ones because the meromorphic dependence on $z$ of the state is \emph{not} inherited by real-quadratic outputs. Indeed, complex conjugation is not an analytic operation. For instance, consider the functional in \cref{eq:bfunctionalq}, applied to a solution map $u$ of the form \cref{eq:frm}. We obtain
\begin{equation}\label{eq:meromorphicsquareexact}
y(z)=\int_\omega|u(z)|^2=\sum_{j,j'=1}^\infty\sum_{i,i'=1}^{\mu_j,\mu_{j'}} \frac{\int_\omega\zeta_{j,i}\overline{\zeta_{j',i'}}}{(\lambda_j-z)^i\overline{(\lambda_{j'}-z)^{i'}}}.
\end{equation}
The quantity of interest $y(z)$ behaves like a meromorphic function \emph{for real wavenumbers} only in special cases, namely, if the resonances $\lambda_j$ are real. Still, even in such situations, it is not meromorphic in the usual sense, since there is no complex neighborhood of the real axis where it is meromorphic.

In the more general case, we may ask whether it is possible to find a surrogate for the exact quadratic functional \cref{eq:meromorphicsquareexact} using any of the methods proposed above. The answer is actually positive for $L^2(\omega')$-SRI (provided $\omega\subset\omega'$), MRI, and (in many cases) POD. For the sake of simplicity, in the following paragraphs we restrict our attention to $L^2(\omega')$-SRI, with $\omega\subset\omega'$, but our discussion applies also to MRI.

Let \cref{eq:meromorphicsquareexact} be the quantity of interest, and let \cref{eq:traceinterpolatorsurrogate} be the $L^2(\omega')$-SRI surrogate of $v_h=u_h|_{\omega'}$. Then, we observe that, by sesquilinearity,
\begin{align*}
    \int_\omega|u(z)|^2\approx&\int_\omega|u_h(z)|^2=\int_\omega|v_h(z)|^2\approx\int_\omega|v_{h,[N]}^{L^2(\omega')\textup{-SRI}}(z)|^2\\
    =&\left(\sum_{i,i'=0}^N\frac{\sum_{j,j'=1}^S\mathring{P}_{ji}\overline{\mathring{P}}_{j'i'}\int_\omega v_h(z_j)\overline{v_h(z_{j'})}}{(z-\zeta_i)\overline{(z-\zeta_{i'})}}\right)\bigg/\left|\sum_{i=0}^N\frac{q_i}{z-\zeta_i}\right|^2\\
    =&\left(\sum_{i,i'=0}^N\frac{\sum_{j,j'=1}^S\mathring{P}_{ji}\overline{\mathring{P}}_{j'i'}\big(G_h^{[y]}\big)_{j'j}}{(z-\zeta_i)\overline{(z-\zeta_{i'})}}\right)\bigg/\left|\sum_{i=0}^N\frac{q_i}{z-\zeta_i}\right|^2,
\end{align*}
where we have defined the $y$-representative Gramian $G_h^{[y]}$ as
\begin{equation*}
    G_h^{[y]}=
    \begin{bmatrix}
    \int_\omega\left|u_h(z_1)\right|^2 & \cdots & \int_\omega u_h(z_S)\overline{u_h(z_1)} \\
    \vdots & \ddots & \vdots\\
    \int_\omega u_h(z_1)\overline{u_h(z_S)} & \cdots & \int_\omega\left|u_h(z_S)\right|^2 
    \end{bmatrix}
    \in\C^{S\times S}.
\end{equation*}
In practical terms, this means that an online-efficient surrogate for $y$ can be built from the $L^2(\omega')$-SRI one, as long as the $y$-representative Gramian is assembled offline. Such assembly can be carried out by exploiting the strategy described in \cref{sec:gramianoverlayed}. Notably, we remark that $G_h^{[y]}$ can be built in conjunction to the $L^2(\omega')$-Gramian $G_h^{(v)}$ in \cref{eq:interpolationsnapgram}: in this case, the same overlay of two meshes can be used to compute $(G_h^{(v)})_{jj'}$, $(G_h^{(v)})_{j'j}$, $(G_h^{[y]})_{jj'}$, and  $(G_h^{[y]})_{j'j}$, for $1\leq j<j'\leq S$.

Note that other similar approaches for handling real-quadratic outputs have been considered in ``system theory'' literature, e.g., (for projection-based MOR) in \cite{quad}.

\section{Numerical results}\label{sec:numerical}

In this section, we showcase the described methods in three numerical examples, where features and limitations of the different techniques emerge. We ran our code in MATLAB\textsuperscript{\textregistered} R2019b on a desktop computer with an 8-core 3.60 GHz Intel\textsuperscript{\textregistered} processor.
The computation of the snapshots via $h$-adaptive FEM is performed using \texttt{p1afem} \cite{FPW11}. The assembly of stiffness and mass matrices associated with pairs of different meshes is carried out using~\cite[Section~6.2]{BPR21}, with a computational complexity that is quadratic in the number of elements of the meshes. 
However, since NVB is a binary refinement rule, it should be possible to design an algorithm that performs this task in at most log-linear complexity with respect to the number of elements.
This artificially inflates the computational times displayed below for the MRI and POD methods.
Numerical testing with more efficient versions of the algorithm is planned for the upcoming future.

\subsection{Toy example on a triangle}

We consider a triangular domain $\Omega=\left\{x\in\R^2,0<x_2<x_1<\frac{\pi}2\right\}$, whose boundary is partitioned into $\Gamma_1\cup\Gamma_2\cup\Gamma_3$. We are interested in the solution $u=u(z)\in H_{\Gamma_1}^1(\Omega)=\{v\in H^1(\Omega),v|_{\Gamma_1}=0\}$ of the Helmholtz equation with uniform forcing term
\begin{equation}
\label{eq:ex1eqn}
\begin{cases}
    -\Delta u(z)-zu(z)=1, &\text{ in } \Omega,\\
    u(z)=0, &\text{ on } \Gamma_1=\left]0,\frac{\pi}2\right[\times\{0\},\\
    \partial_{\bm{\nu}} u(z)=\partial_{x_1} u(z)=0, &\text{ on } \Gamma_2=\{\frac{\pi}2\}\times\left]0,\frac{\pi}2\right[,\\
    \partial_{\bm{\nu}} u(z)=\frac1{\sqrt{2}}\left(\partial_{x_2} u(z)-\partial_{x_1} u(z)\right)=0, &\text{ on } \Gamma_3=\partial \Omega\setminus\{\Gamma_1\cup\Gamma_2\}.
\end{cases}
\end{equation}

The eigenproblem for the minus Laplacian over $\Omega$ with the above boundary conditions can be solved explicitly. The eigenvalues are
\begin{equation}
    \label{eq:ex1eigval}
    \lambda_{m,n}=m^2+n^2,\quad m,n\in\mathbb{N}\text{ odd},m\leq n,
\end{equation}
and the corresponding $L^2(\Omega)$-orthonormal eigenfunctions are
\begin{equation}
    \label{eq:ex1eigfunc}
    \phi_{m,n}(x)=
    \begin{cases}
        \frac{4\sqrt{2}}\pi\sin(mx_1)\sin(nx_2), &\text{ if } m=n,\\
        \frac{4}\pi\left(\sin(mx_1)\sin(nx_2)+\sin(nx_1)\sin(mx_2)\right), &\text{ if } m<n.
    \end{cases}
\end{equation}
Accordingly, the solution $u$ admits the closed-form infinite series expression
\begin{align}
    u(z)|_{x}=&\sum_{m=1,3,\ldots}\frac{2\sqrt{2}\phi_{m,m}(x)}{\pi m^2(\lambda_{m,m}-z)}+\sum_{\substack{m,n=1,3,\ldots\\m<n}}\frac{4\phi_{m,n}(x)}{\pi mn(\lambda_{m,n}-z)}\nonumber\\
    =&\sum_{m,n=1,3,\ldots}\frac{16}{\pi^2mn(m^2+n^2-z)}\sin(mx_1)\sin(nx_2),
    \label{eq:ex1exact}
\end{align}
where the series coefficients have been obtained by $L^2(\Omega)$-projection of the forcing term onto the eigenspaces. Hence, in the scope of evaluating the FEM and MOR errors in post-processing, the reference solution $u(z)$ is available to arbitrary precision by truncation of \cref{eq:ex1exact}.

\subsubsection{Adaptive FEM approximation for fixed frequency}\label{sec:numhfem}

In our first experiment,
we present the computation of a snapshot
using the $h$-adaptive FEM described in \cref{sec:hfem}.
To this end,
we fix $z=k^2=51$ 
and run \cref{alg:hfem}
with $\theta=0.1$, $\tol_h=5\cdot 10^{-2}$, and 
$\Nmax =\frac{\abs{\Omega}k^4}{4\tol_h^2}$.
With this choice, we are heuristically trying to guarantee that the finest possible mesh
used by the algorithm satisfies a ``resolution condition'' of the form $\widetilde{h}_\bullet k^2\sim 2\tol_h$ \cite{bisg1997a},
where $\widetilde{h}_\bullet^2=\abs{\Omega}/N_\bullet$ is the mesh size of a 2D uniform mesh made of right triangles.
In this specific case, we have $\Nmax \approx 3.2\cdot 10^6$, but the adaptive algorithm stops beforehand at $N_{143}\approx 1.9\cdot 10^5$, since $\tol_h$ is attained at the 143\textsuperscript{rd} iteration.
Note that $z=51$ is not of the form \cref{eq:ex1eigval}, and its closest resonance is $\lambda_{1,7}=\lambda_{5,5}=50$.

\begin{figure}[h!]
\centering
    \begin{minipage}[t]{.4\textwidth}
     \begin{subfigure}[t]{.925\textwidth}
         \includegraphics[width=\textwidth,trim=1cm .65cm 1cm .75cm, clip]{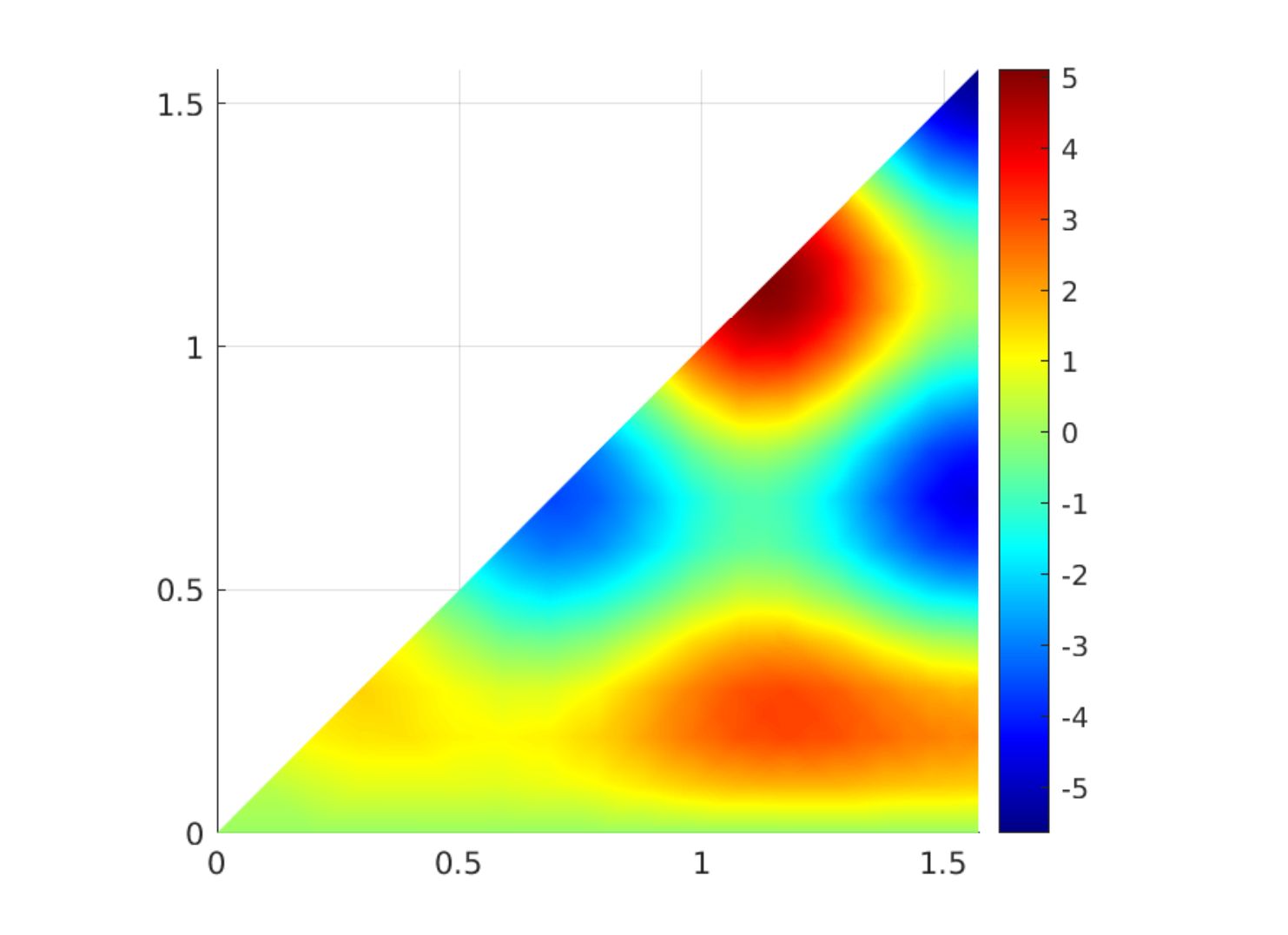}
         \caption{$u_{50}^\textup{FEM}(z)$}
         \label{fig:ex1f1s051}
     \end{subfigure}
     
     \vspace{.5mm}
     \begin{subfigure}[t]{\textwidth}
         \pgfplotstableread[col sep=comma]{Data/ex1estimator_long.csv}\tikzdataires
         \hspace{-.1cm}%
         \begin{tikzpicture}
         \begin{loglogaxis}[
            ticklabel style = {font=\tiny},
            width = .95\textwidth, height = .8\textwidth,
            xlabel={\small DoFs $N_\bullet$},
            xlabel style={yshift=.1cm},
            legend style={nodes={scale=0.75, transform shape}, at={(.75,.975)}},
            xmin = 1e0, xmax = 9e6,
            ymin = 1e-3, ymax = 1e3
            ]
            \addplot[black, mark=*, mark size=.1] plot table [x=dofs, y=est]{\tikzdataires};
            \addplot[blue, semithick, densely dotted] plot table [x=dofs, y=err]{\tikzdataires};
            \draw[black, densely dashed] (axis cs:1e4,1)--(axis cs:4e6,.05);
            \draw[red, densely dotted] (axis cs:1e-1,5e-2)--(axis cs:1e8,5e-2);
            \draw[blue, densely dotted] (axis cs:3.2e6,1e-4)--(axis cs:3.2e6,1e4);
            \draw (axis cs:2.5e5,.7) node[scale=0.65, rotate=-30] {$\mathcal{O}(N_\bullet^{-1/2})$};
            \draw[red] (axis cs:10,2e-2) node[scale=0.85] {$\tol_h$};
            \draw[blue] (axis cs:1e6,1e2) node[scale=0.85, rotate=90] {$N_\text{max}$};
            \legend{$\eta_\bullet(\Omega)$,$e_\bullet(\Omega)$}
         \end{loglogaxis}
         \end{tikzpicture}
         \vspace{-.125cm}
         \caption{$\eta_\bullet(\Omega)$ and $e_\bullet(\Omega)$ vs.\ FEM DoFs.}
         \label{fig:ex1f1conv}
     \end{subfigure}
    \end{minipage}%
    \begin{minipage}[t]{.4\textwidth}
     \begin{subfigure}[t]{.925\textwidth}
         \includegraphics[width=\textwidth,trim=1cm .65cm 1cm .75cm, clip]{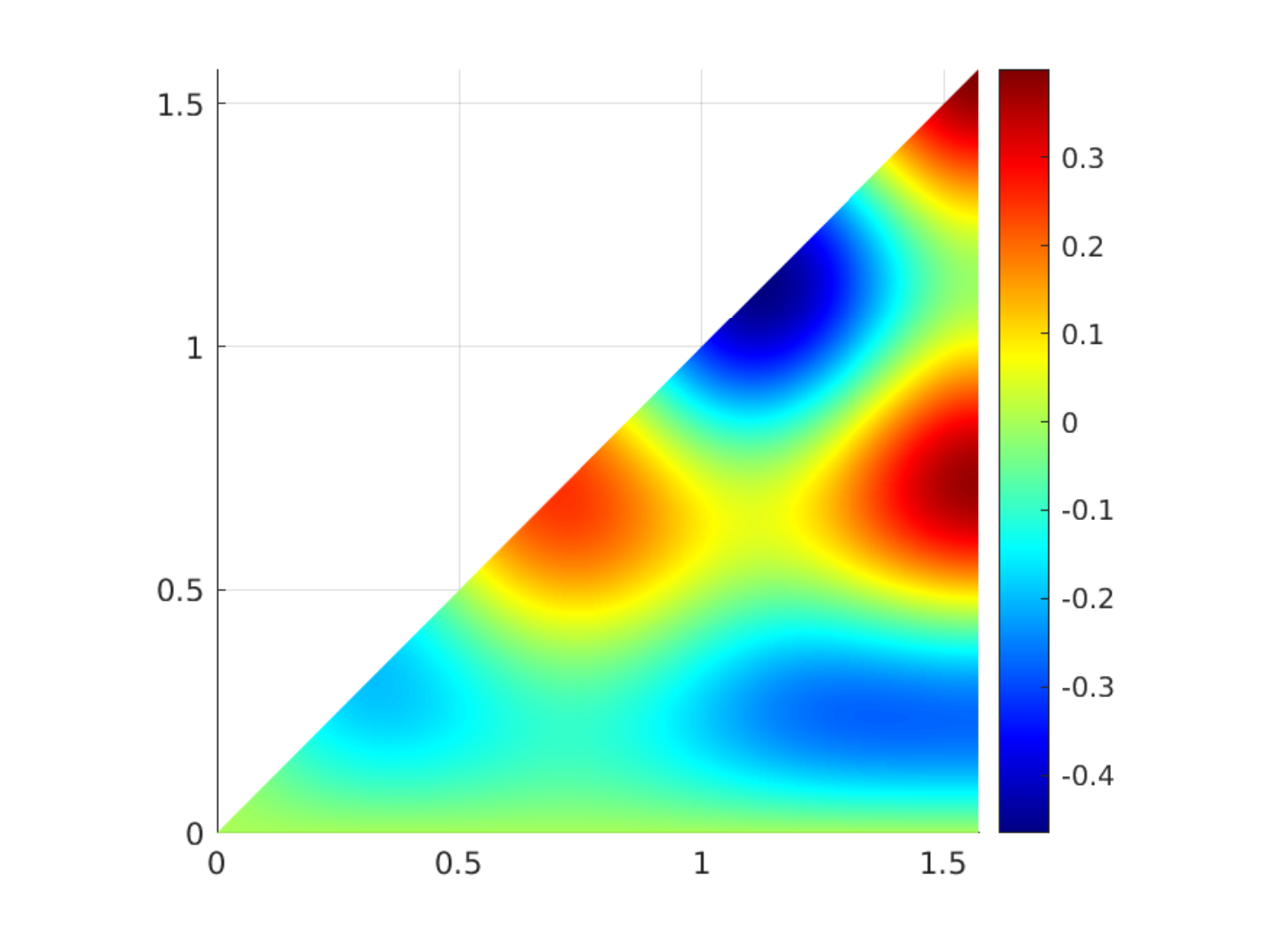}
         \caption{$u_{143}^\textup{FEM}(z)$}
         \label{fig:ex1f1s143}
     \end{subfigure}
          
     \begin{subfigure}[t]{.925\textwidth}
         \includegraphics[width=\textwidth,trim=1cm .65cm 1cm .75cm, clip]{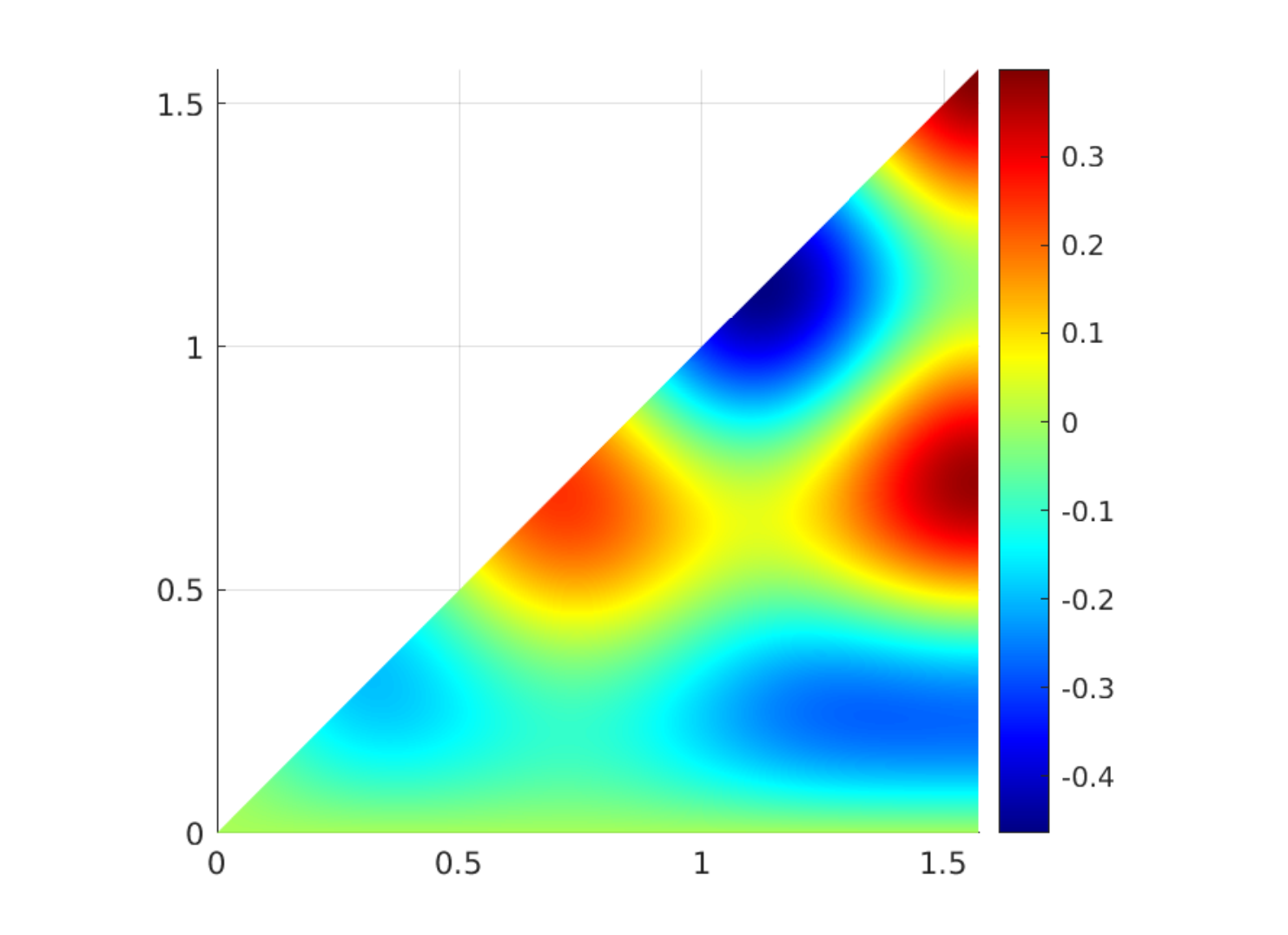}
         \caption{Exact $u(z)$}
         \label{fig:ex1f1exact}
     \end{subfigure}
    \end{minipage}
     \caption{$h$-adaptive FEM results for $z=51$.}
     \label{fig:ex1f1}

    \begin{minipage}[t]{.45\textwidth}
        \vspace{0pt}%
         \pgfplotstableread[col sep=comma]{Data/ex1trace.csv}\tikzdataitr
         \begin{tikzpicture}
         \begin{axis}[
            ticklabel style = {font=\tiny},
            width = \textwidth, height = .575\textwidth,
            xlabel={\small $x_2$}, ylabel={\small $v(51)|_{x_2}$},
            xlabel style={yshift=.1cm}, ylabel style={yshift=-.25cm},
            legend columns=3,
            ymin = -1, ymax = 2
            ]
            \addplot[black] plot table [y=exact]{\tikzdataitr};
            \addplot[blue, densely dotted] plot table [y=coarse]{\tikzdataitr};
            \addplot[red, densely dashed] plot table [y=fine]{\tikzdataitr};
            \legend{exact,$v_{50}$,$v_{143}$};
         \end{axis}
         \end{tikzpicture}
         \caption{$u(51)|_{(\frac\pi2,x_2)}$.}
         \label{fig:ex1trace}
    \end{minipage}%
    \begin{minipage}[t]{.5\textwidth}
        \vspace{.5cm}%
        \centering
        \bgroup
        \setlength{\tabcolsep}{2pt}
        \begin{tabular}{c|c|c|c|}
            & exact & rel. err. & rel. err. \\
            & value & on $\mathcal{T}_{50}$ & on $\mathcal{T}_{143}$ \\[.1cm]
            \hline
            $y(51)$ & \texttt{1.47e-1} & \texttt{1.76e+1} & \texttt{3.19e-3} \\
            \hline
        \end{tabular}
        \egroup
        \vspace{.325cm}
        \captionof{table}{QoI and relative FEM errors.}
        \label{tab:ex1qois}
    \end{minipage}
\end{figure}

In~\cref{fig:ex1f1},
we compare the analytical solution (\cref{fig:ex1f1exact}) with the FE approximations
at the 50\textsuperscript{th} (\cref{fig:ex1f1s051}) and at the 143\textsuperscript{rd} iteration (\cref{fig:ex1f1s143}).
In \cref{fig:ex1f1conv}, we show the evolution of the error estimator $\eta_\bullet$ and the true error $e_\bullet(z)=\norm{\nabla(u_\bullet(z)-u(z))}{L^2(\Omega)}$ as the mesh gets adaptively refined.
Note that, in \cref{alg:hfem}, the $h$-adaptive loop stops as soon as $\eta_\bullet$ becomes smaller than the tolerance $\tol_h$. In \cref{fig:ex1f1conv}, we are also including points that correspond to further
adaptive steps of refinement only for illustrative reasons.

Several peaks appear before the asymptotic (optimal) convergence regime is reached. Such peaks are caused by resonances of the discrete problem at the chosen value of $z$. Indeed, since $z$ is \emph{not} a resonance of the continuous problem \cref{eq:ex1eqn}, we know that the discrete problem will not have a resonance at $z$ \emph{for a fine enough mesh}. Still, resonances of the discrete problem may occur at $z$ if the mesh is too coarse. For accuracy, it is crucial that the adaptive algorithm stops after all peaks. We showcase this in \cref{fig:ex1f1s051}, where we show the intermediate FEM solution after $50$ iterations, i.e., just to the left of the final peak in \cref{fig:ex1f1conv}. We observe that the approximation of $u$ is rather poor, because the discrete resonance closest to $z$ is placed \emph{on the wrong side of $z$}.

Instead, when the convergence enters the asymptotic regime, the error estimator and the true error decay with the same (optimal) rate. This is in agreement with the equivalence established in \cref{eq:Ceff} (note that here the oscillation terms vanish, because $f$ and $g_N$ are constant in this example).

Given $u$ and its FE approximation, we introduce the QoI
\begin{equation}
\label{eq:ex1qoi1}
    y(z) = \int_{\Gamma_2}u(z)=\left(u(z),1\right)_{L^2(\Gamma_2)}. 
\end{equation}
and the restriction $v(z):= u(z)\vert_{\Gamma_2}$ to $\Gamma_2$ (the domain of integration in~\cref{eq:ex1qoi1}).
In \cref{fig:ex1trace},
we compare $v(z)$ with its FE approximation obtained by restricting $u^\textup{FEM}_{51}(z)$ and $u^\textup{FEM}_{143}(z)$ to $\Gamma_2$.
We observe a complete mismatch between the exact solution and the coarse approximation,
confirming our qualitative observations from the previous paragraphs.
In~\cref{tab:ex1qois}, we make our conclusions quantitative by looking at the relative approximation error in the approximation of the QoI.

\subsubsection{Surrogate model for a linear quantity of interest}

Now, we move to the approximation of the linear QoI $y(z)$ in \cref{eq:ex1qoi1} over the range $z\in Z=[1,100]$.
First, we build a rational surrogate of type $[14/14]$ by SRI, using snapshots at $S=29$ uniformly spaced sample points in $Z$. Specifically, we compute the snapshots by using $h$-adaptive FEM with the same parameters ($\theta$, $\tol_h$, and $N_\text{max}$) as in the previous section. Note, in particular, that $N_\text{max}$ increases with $z$. When computing the snapshots at $z\approx 25.75$, we run into an issue: the $h$-adaptive loop in \cref{alg:hfem} is terminated before the prescribed $\tol_h$ is attained, because the maximum number of elements is reached. This is caused by the nearby presence of the eigenvalue $\lambda_{1,5}=26$, which slows down the convergence of the FEM error. At this point, we are faced with a choice:
\begin{itemize}
    \item We can use the inaccurate non-converged snapshot as it is. This is the simplest option, but might result in a poor surrogate, in particular near the affected sample point.
    \item We can discard the non-converged snapshot and build the surrogate without it. This is the ``wasteful'' option, since we are ignoring the results of valuable offline computations.
    \item We can (temporarily) increase $N_\text{max}$, e.g., multiplying its value by 10, and restart the $h$-adaptive iterations in the hope of convergence. If the increased $N_\text{max}$ leads to convergence, we use the snapshot. Otherwise, we discard it. This is the most robust, but potentially expensive, option. Here, we follow this approach.
\end{itemize}
Note that it is not necessarily a good idea to ditch $N_\text{max}$ altogether, or, equivalently, to set its value to $\infty$. Indeed, due to the \emph{a priori} unknown location of the poles, it might happen that a sample point is extremely close to a pole of $u$. Correspondingly, the $h$-adaptive loop will necessarily take a long time to converge there. As such, we believe it more appropriate to be conservative, by setting a reasonably large (but finite) $N_\text{max}$, if necessary retrying with $10N_\text{max}$, and throwing away the snapshot in case of non-convergence.

In addition to the type $[14/14]$ SRI, we also build a type $[14/14]$ MRI and a type $[15/15]$ POD surrogate. Building the latter two reduced models requires only 15 snapshots, as opposed to the 29 needed for SRI. We take such snapshots at 15 uniformly spaced values of $z\in Z$. Note that, obviously, the issue of non-converging snapshots can appear also in MRI and POD. We deal with it as described above.

\begin{figure}[t!]
     \pgfplotstableread[col sep=comma]{Data/ex2global.csv}\tikzdatag
     \pgfplotstableread[col sep=comma]{Data/ex2local.csv}\tikzdatal
     \centering
     \begin{tikzpicture}
     \begin{semilogyaxis}[
        scale only axis, grid,
        ticklabel style = {font=\tiny},
        width = .5\textwidth, height = 2cm,
        ylabel={\small $\abs{y(z)}$},
        ylabel style={yshift=-.25cm},
        xticklabels = {,,},
        legend style={nodes={scale=0.6, transform shape}},
        legend columns=2,
        xmin = 1, xmax = 100,
        name = top left plot
        ]
        \addplot[blue] plot table [x=z, y=sri]{\tikzdatag};
        \addplot[red, densely dashed] plot table [x=z, y=mri]{\tikzdatag};
        \addplot[green!70!black, densely dotted] plot table [x=z, y=pod]{\tikzdatag};
        \addplot[black, only marks, mark = *, mark size = .1mm] plot table [x=z, y=ex]{\tikzdatag};
        \legend{SRI,MRI,POD,analytic};
     \end{semilogyaxis}

     \begin{semilogyaxis}[
        scale only axis, grid,
        ticklabel style = {font=\tiny},
        width = .5\textwidth, height = 2.375cm,
        xlabel={\small $z$}, ylabel={\scriptsize $\abs{\widetilde{y}(z)-y(z)}/\abs{y(z)}$},
        xlabel style={yshift=.1cm}, ylabel style={yshift=-.25cm},
        xmin = 1, xmax = 100,
        ytick = {1e-6,1e-4,1e-2,1e0,1e2},
        ymax = 7e1,
    	at={(top left plot.below south west)},
    	anchor=north west
        ]
        \addplot[blue] plot table [x=z, y=sri_eR]{\tikzdatag};
        \addplot[red, densely dashed] plot table [x=z, y=mri_eR]{\tikzdatag};
        \addplot[green!70!black, densely dotted] plot table [x=z, y=pod_eR]{\tikzdatag};
        \addplot[black, only marks, mark = +, mark size = .75mm, domain = 1:99.9, samples = 29] {3e1};
        \addplot[black, only marks, mark = x, mark size = .75mm, domain = 1:99.9, samples = 15] {3e1};
     \end{semilogyaxis}

     \begin{semilogyaxis}[
        scale only axis, grid,
        ticklabel style = {font=\tiny},
        width = .25\textwidth, height = 2cm,
        xticklabels = {,,},
        legend style={nodes={scale=0.6, transform shape}},
        xmin = 33, xmax = 35,
        ymin = 5e-4,
    	at={(top left plot.right of north east)},
    	anchor=left of north west,
    	name = top right plot
        ]
        \addplot[blue] plot table [x=z, y=sri]{\tikzdatal};
        \addplot[red, densely dashed] plot table [x=z, y=mri]{\tikzdatal};
        \addplot[green!70!black, densely dotted] plot table [x=z, y=pod]{\tikzdatal};
        \addplot[black, only marks, mark = *, mark size = .1mm] plot table [x=z, y=ex]{\tikzdatal};
        \legend{SRI,MRI,POD,analytic};
     \end{semilogyaxis}

     \begin{semilogyaxis}[
        scale only axis, grid,
        ticklabel style = {font=\tiny},
        width = .25\textwidth, height = 2.375cm,
        xlabel={\small $z$},
        xlabel style={yshift=.1cm},
        xmin = 33, xmax = 35,
        ytick = {1e-6,1e-4,1e-2,1e0,1e2},
    	at={(top right plot.below south west)},
    	anchor=north west
        ]
        \addplot[blue] plot table [x=z, y=sri_eR]{\tikzdatal};
        \addplot[red, densely dashed] plot table [x=z, y=mri_eR]{\tikzdatal};
        \addplot[green!70!black, densely dotted] plot table [x=z, y=pod_eR]{\tikzdatal};
     \end{semilogyaxis}
     \end{tikzpicture}
    \caption{Surrogate $\abs{y}$ with analytic validation points (top row) and relative approximation error (bottom row). The plots in the right column are a zoom on $z\in[33,35]$. On the top of the bottom left plot, we show with stars and pluses the locations of the sample points (pluses are exclusive to SRI).\vspace{-2mm}}
    \label{fig:ex1f2}
\end{figure}
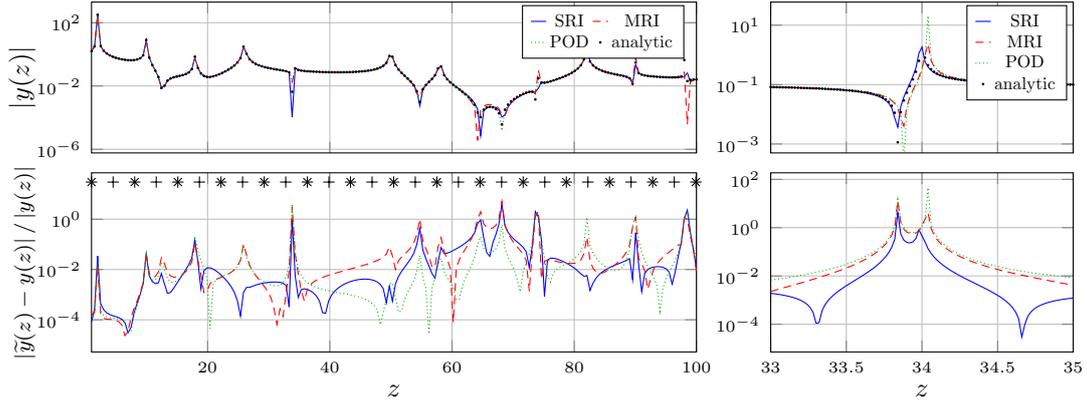

We show the results of the approximation in \cref{fig:ex1f2}. Note that the approximation error is computed with respect to the analytic solution \cref{eq:ex1exact}, so that, in particular, the error does not vanish at the sample points even though all approaches are interpolatory. This is due to the baseline FEM error, which acts as noise. We see that the approximations yielded by the three approaches are quite similar, with MRI behaving only slightly worse (on average) than the other two methods.

Looking locally around the (arbitrarily chosen) pole at $\lambda_{3,5}=34$, we see that SRI approximates best the location of the pole. This should be expected, due to the higher sampling resolution. Indeed, locally around $\lambda_{3,5}$, SRI uses samples at both $z\approx 32.82$ and $z\approx 36.36$ while MRI and POD only have the one at $z\approx 36.36$.

In \cref{tab:ex1time} we show the offline time of the three approaches. The computation of the snapshots in SRI takes considerably longer than for MRI/POD due to the non-converging samples discussed above. More generally, if the $S$ MRI/RB snapshots are a subset of the $(2S-1)$ SRI ones, the sampling for SRI will obviously take more time, but not necessarily twice as much. Overall, we see that building the SRI and MRI surrogates is about 20\% faster than POD. Note that, if the SRI samples had all converged at the first try, we could have expected the offline phase of SRI to be, overall, about 4 times faster than that of MRI.

\begin{figure}[h!]
\centering
\begin{small}
\begin{tabular}{c|c|c|c|}
Method & SRI & MRI & POD\\
\hline
\hline
Snapshots & \texttt{8.21e+2} & \multicolumn{2}{c|}{\texttt{1.06e+2}}\\
\hline
Build & & & \\
snapshot & \texttt{-------} & \texttt{7.00e+2} & \texttt{1.04e+3}\\
Gramian(s) & & & \\
\hline
Assemble & \multirow{2}{*}{\texttt{3.90e-2}} & \multirow{2}{*}{\texttt{1.81e-2}} & \multirow{2}{*}{\texttt{1.60e-2}}\\
surrogate & & & \\
\hline
\hline
Total & \texttt{8.21e+2} & \texttt{8.07e+2} & \texttt{1.15e+3}\\
\hline
\end{tabular}
\end{small}
\captionof{table}{Timings in seconds for the construction of $\widetilde{y}$.\vspace{-2mm}}
\label{tab:ex1time}
\end{figure}

\subsection{Vibrations of a non-convex elastic plate}

\begin{figure}[p]
\centering
    \begin{minipage}[b]{.225\textwidth}
         \hspace{0.cm}
         \begin{flushleft}
         \includegraphics[width = .9\textwidth,trim=1cm .65cm 1cm .75cm, clip]{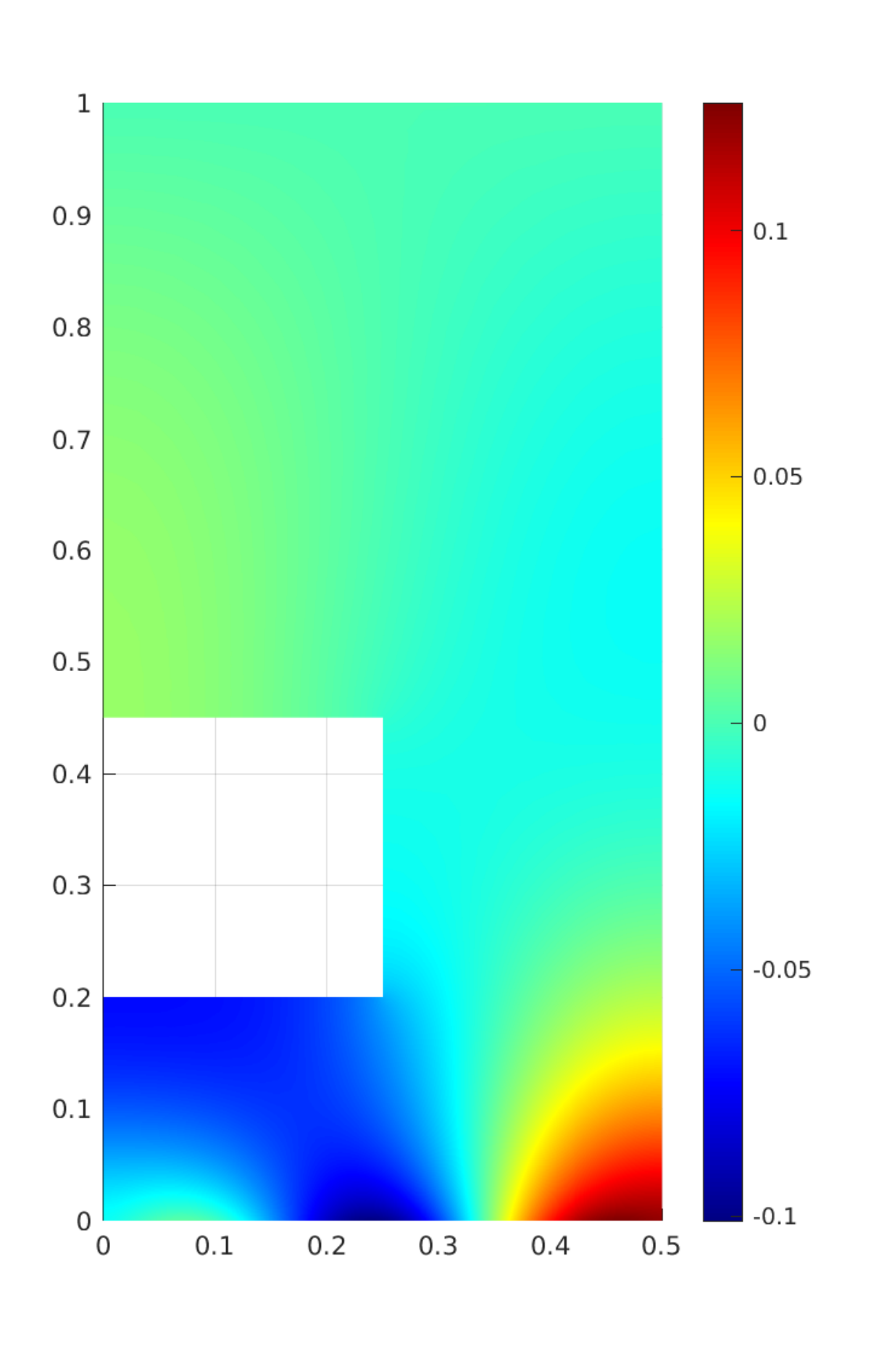}
         \end{flushleft}
     \begin{subfigure}[t]{.975\textwidth}
         \includegraphics[width = .785\textwidth,trim=1cm .65cm 1cm .75cm, clip]{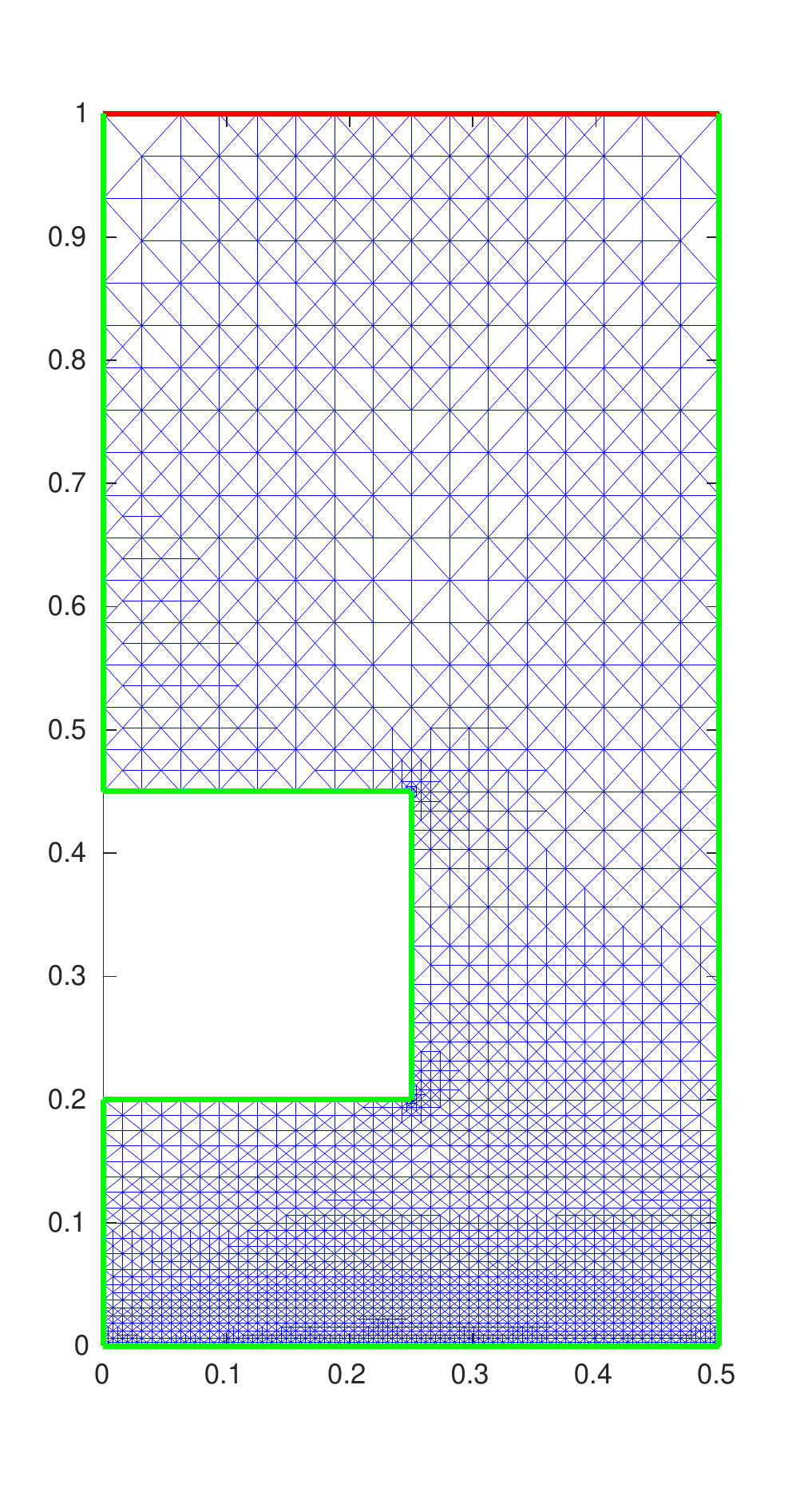}
         \caption{$u_{70}^\textup{FEM}(50)$}
         \label{fig:exHm}
    \end{subfigure}
    \end{minipage}%
    \begin{minipage}[b]{.0375\textwidth}
        \phantom{lll}
    \end{minipage}%
    \begin{minipage}[b]{.225\textwidth}
         \hspace{0.cm}
         \begin{flushleft}
         \includegraphics[width = .9\textwidth,trim=1cm .65cm 1cm .75cm, clip]{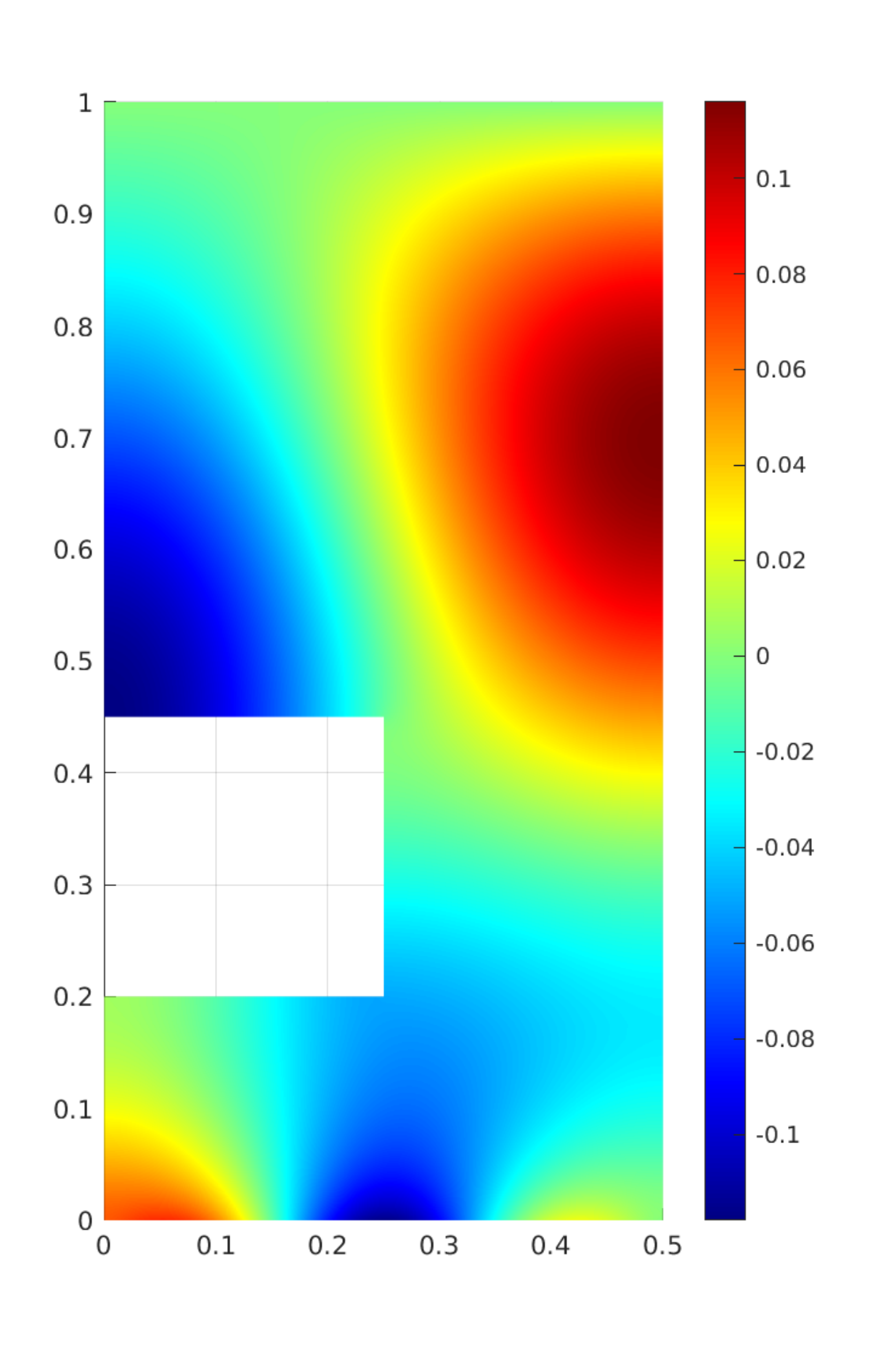}
         \end{flushleft}
     \begin{subfigure}[t]{.975\textwidth}
         \includegraphics[width = .785\textwidth,trim=1cm .65cm 1cm .75cm, clip]{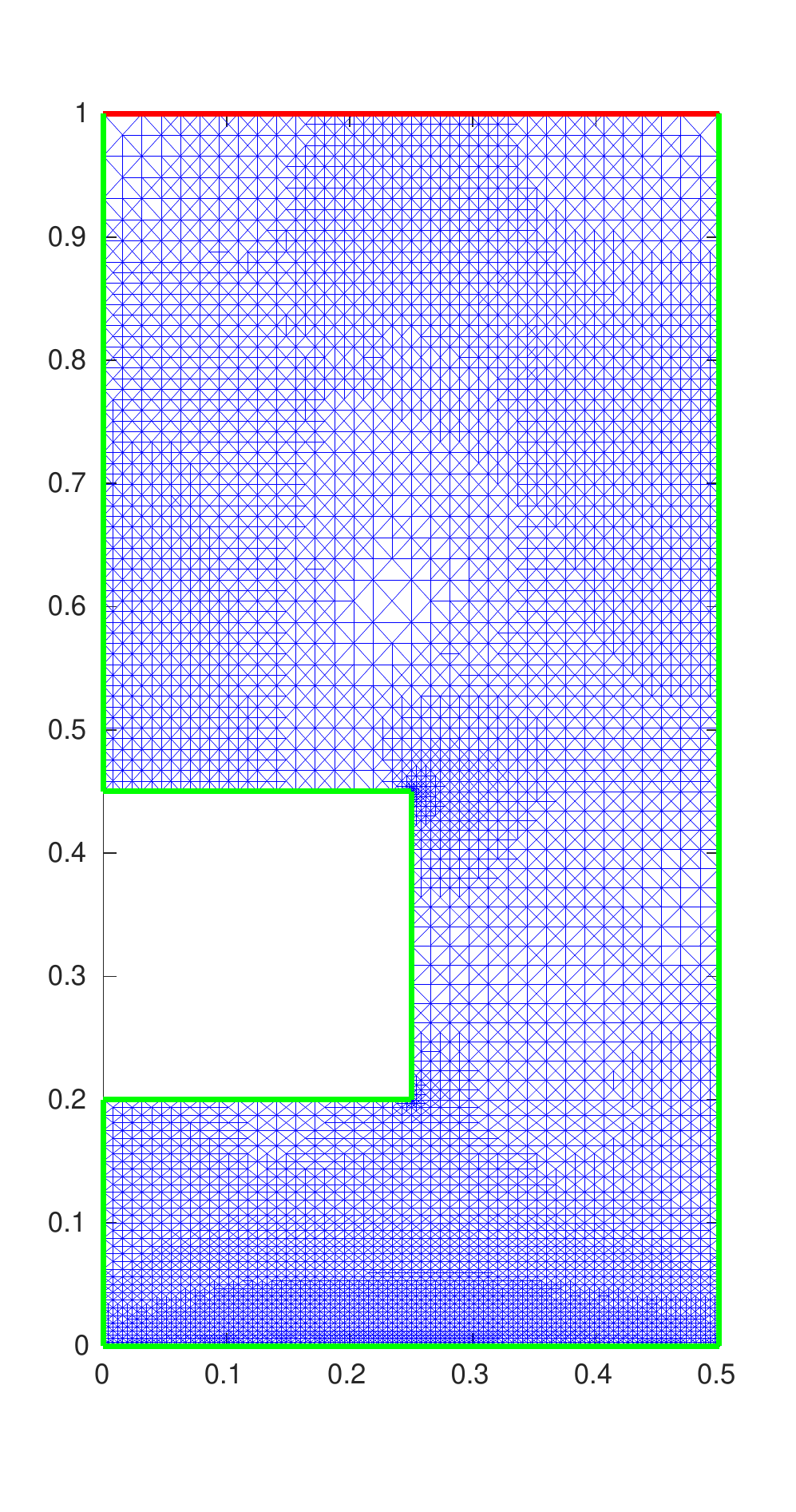}
         \caption{$u_{90}^\textup{FEM}(55)$}
         \label{fig:exHl}
    \end{subfigure}
    \end{minipage}%
    \begin{minipage}[b]{.475\textwidth}
         \hspace{0.cm}
     \begin{subfigure}[t]{.975\textwidth}
         \pgfplotstableread[col sep=comma]{Data/ex5estimator.csv}\tikzdataires
         \centering
         \begin{tikzpicture}
         \begin{loglogaxis}[
            ticklabel style = {font=\tiny},
            width = .9\textwidth, height = \textwidth,
            xlabel={\small DoFs $N_\bullet$}, ylabel={\small $\eta_\bullet(\Omega)$},
            xlabel style={yshift=.1cm}, ylabel style={yshift=-.25cm},
            ymin = 4e-2, ymax = 2e1
            ]
            \addplot[semithick, black, mark=*, mark size=.2] plot table [x=dofs50, y=est50]{\tikzdataires};
            \addplot[semithick, blue, densely dotted] plot table [x=dofs55, y=est55]{\tikzdataires};
            \draw[black, densely dashed] (axis cs:1e3,.3)--(axis cs:16e3,.075);
            \draw (axis cs:6e3,.2) node[scale=.85, rotate=-30] {$\mathcal{O}(N_\bullet^{-1/2})$};
            \legend{$z=50$,$z=55$};
         \end{loglogaxis}
         \end{tikzpicture}
         \caption{$\eta_\bullet(\Omega)$ vs.\ FEM DoFs.}
         \label{fig:exHr}
     \end{subfigure}
     \vspace{10pt}
     
     \begin{subfigure}[t]{.975\textwidth}
         \pgfplotstableread[col sep=comma]{Data/ex5trace.csv}\tikzdataitr
         \centering
         \begin{tikzpicture}
         \begin{axis}[
            ticklabel style = {font=\tiny},
            width = .925\textwidth, height = .625\textwidth,
            xlabel={\small $x_1$}, ylabel={\small $v(z)|_{x_1}$},
            xlabel style={yshift=.1cm}, ylabel style={yshift=-.25cm},
            legend style={nodes={scale=0.6, transform shape}},
            legend pos=north west,
            legend columns=2,
            xmin = 0, xmax = .5,
            ymin = -.15, ymax = .175
            ]
            \addplot[black] plot table [y=tr50]{\tikzdataitr};
            \addplot[blue, densely dotted] plot table [y=tr55]{\tikzdataitr};
            \legend{$v_{70}(50)$,$v_{90}(55)$};
         \end{axis}
         \end{tikzpicture}
         \caption{$u(z)|_{(x_1,1)}$.}
         \label{fig:exHtrace}
     \end{subfigure}
    \end{minipage}
     \caption{$h$-adaptive FEM results for $z=50$ and $z=55$.}
     \label{fig:exHs}
\end{figure}

We consider a rectangular domain with a square hole 
\begin{equation*}
\Omega=\left]0,0.5\right[\times\left]0,1\right[\setminus\left]0,0.25\right[\times\left]0.2,0.45\right[.
\end{equation*}
We denote by $\Gamma_1=\left]0,0.5\right[\times\{0\}$ and $\Gamma_2=\left]0,0.5\right[\times\{1\}$ the bottom and top sides of $\Omega$, respectively, see \cref{fig:exHs}. We are interested in the solution $u=u(z)\in H_{\Gamma_2}^1(\Omega)=\{v\in H^1(\Omega),v|_{\Gamma_2}=0\}$ of the Helmholtz equation
\begin{equation}
\label{eq:exHeqn}
\begin{cases}
    -\Delta u(z)-zu(z)=0, &\text{ in } \Omega,\\
    \partial_{\bm{\nu}} u(z)=g(z), &\text{ on } \Gamma_1,\\
    u(z)=0, &\text{ on } \Gamma_2,\\
    \partial_{\bm{\nu}} u(z)=0, &\text{ on } \partial\Omega\setminus\{\Gamma_1\cup\Gamma_2\}.
\end{cases}
\end{equation}

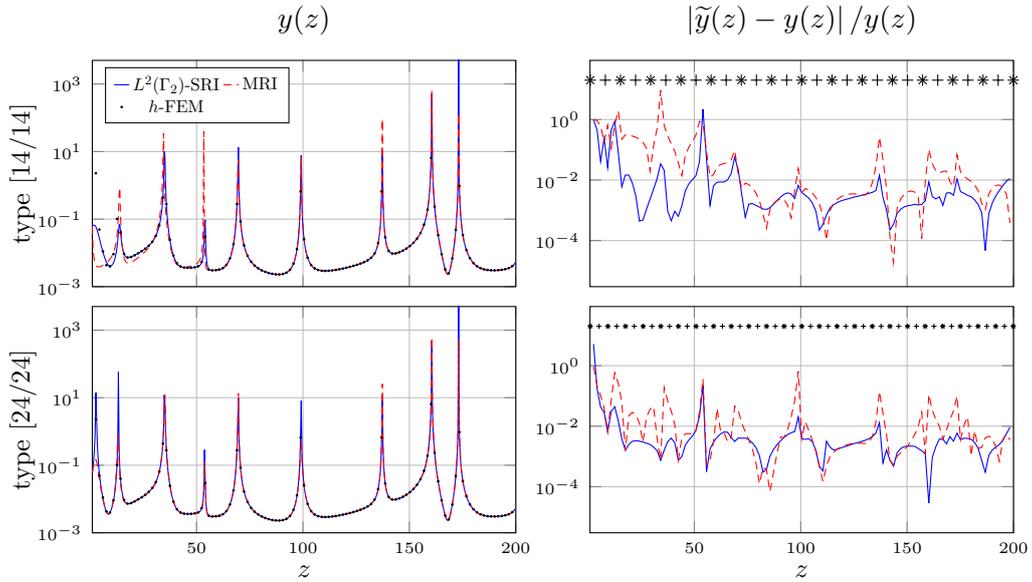
\begin{figure}[p]
     \pgfplotstableread[col sep=comma]{Data/ex6global.csv}\tikzdatag
     \centering
     \begin{tikzpicture}
     \begin{semilogyaxis}[
        scale only axis, grid,
        ticklabel style = {font=\tiny},
        width = .35\textwidth, height = 3cm,
        ylabel={\small type $[14/14]$},
        ylabel style={yshift=-.25cm},
        xticklabels = {,,},
        xmin = 1, xmax = 200,
        ytick = {1e-3,1e-1,1e1,1e3},
        ymin = 1e-3, ymax = 5e3,
        legend pos=north west,
        legend columns = 2,
        legend style={nodes={scale=0.6, transform shape}},
        name = top left plot,
        title = {$y(z)$}
        ]
        \addplot[blue] plot table [x=z, y=vri29]{\tikzdatag};
        \addplot[red, densely dashed] plot table [x=z, y=mri15]{\tikzdatag};
        \addplot[black, only marks, mark = *, mark size = .1mm] plot table [x=z, y=ex]{\tikzdatag};
        \legend{$L^2(\Gamma_2)$-SRI,MRI,$h$-FEM};
     \end{semilogyaxis}

     \begin{semilogyaxis}[
        scale only axis, grid,
        ticklabel style = {font=\tiny},
        width = .35\textwidth, height = 3cm,
        xlabel={\small $z$},
        xlabel style={yshift=.1cm},
        ylabel={\small type $[24/24]$},
        ylabel style={yshift=-.25cm},
        xmin = 1, xmax = 200,
        ytick = {1e-3,1e-1,1e1,1e3},
        ymin = 1e-3, ymax = 5e3,
    	at={(top left plot.below south west)},
    	anchor=north west
        ]
        \addplot[blue] plot table [x=z, y=vri49]{\tikzdatag};
        \addplot[red, densely dashed] plot table [x=z, y=mri25]{\tikzdatag};
        \addplot[black, only marks, mark = *, mark size = .1mm] plot table [x=z, y=ex]{\tikzdatag};
     \end{semilogyaxis}

     \begin{semilogyaxis}[
        scale only axis, grid,
        ticklabel style = {font=\tiny},
        width = .35\textwidth, height = 3cm,
        xmin = 1, xmax = 200,
        xticklabels = {,,},
        ytick = {1e-4,1e-2,1e0},
        ymin = 3e-6, ymax = 9e1,
    	at={(top left plot.right of north east)},
    	anchor=left of north west,
        name = top right plot,
        title={$\abs{\widetilde{y}(z)-y(z)}/y(z)$}
        ]
        \addplot[blue] plot table [x=z, y=vri29_eR]{\tikzdatag};
        \addplot[red, densely dashed] plot table [x=z, y=mri15_eR]{\tikzdatag};
        \addplot[black, only marks, mark = +, mark size = .75mm, domain = 1.1:199.9, samples = 29] {2e1};
        \addplot[black, only marks, mark = x, mark size = .75mm, domain = 1.1:199.9, samples = 15] {2e1};
     \end{semilogyaxis}

     \begin{semilogyaxis}[
        scale only axis, grid,
        ticklabel style = {font=\tiny},
        width = .35\textwidth, height = 3cm,
        xlabel={\small $z$},
        xlabel style={yshift=.1cm},
        xmin = 1, xmax = 200,
        ytick = {1e-4,1e-2,1e0},
        ymin = 3e-6, ymax = 9e1,
    	at={(top right plot.below south west)},
    	anchor=north west
        ]
        \addplot[blue] plot table [x=z, y=vri49_eR]{\tikzdatag};
        \addplot[red, densely dashed] plot table [x=z, y=mri25_eR]{\tikzdatag};
        \addplot[black, only marks, mark = +, mark size = .35mm, domain = 1.1:199.9, samples = 49] {2e1};
        \addplot[black, only marks, mark = x, mark size = .35mm, domain = 1.1:199.9, samples = 25] {2e1};
     \end{semilogyaxis}
    \end{tikzpicture}
    \caption{Surrogates for $y$ with high-fidelity $h$-FEM validation points (left column) and relative approximation errors $\abs{\widetilde{y}-y}/y$ (right column). On top of the right plots, we show with stars and pluses the locations of the sample points (pluses are exclusive to $L^2(\Gamma_1)$-SRI). The rational type increases moving down the plots.}
    \label{fig:exHf}
\end{figure}

The solution $u$ corresponds to the transverse displacement field of a thin membrane $\Omega$ in the ``small deformation'' regime, when the membrane, clamped at $\Gamma_2$, is subject to a time-harmonic excitation. The Neumann forcing term
\begin{equation*}
    g(z)|_x=-0.15\text{i}\sqrt{z}\exp\left(\frac{3\sqrt{3}}2\text{i}\sqrt{z}x_1\right)
\end{equation*}
denotes the force exerted on the membrane by the plane wave $u_\text{inc}(z)|_x=0.1\exp\left(3\text{i}\sqrt{z}x\cdot\widehat{x}\right)$ (with $\widehat{x}=(\cos(\frac\pi6),\sin(\frac\pi6))^\top$ being the direction of propagation of the wave) that impinges on $\Gamma_1$ from below. Note that $g$ is not affine in $z$, cf.\ \cref{eq:affine}, preventing a straightforward application of POD. We set the localized real-quadratic QoI
\begin{equation}\label{eq:exHqoi}
    y(z)=\int_{\Gamma_1}|u(z)|^2=\|u(z)\|_{L^2(\Gamma_1)}^2,
\end{equation}
i.e., the mean-square displacement on $\Gamma_1$, as target of the approximation endeavor. This leads us to the definition of the intermediate quantity $v(z)=u(z)|_{\Gamma_1}\in L^2(\Gamma_1)$.

To obtain the FEM solution, we employ the $h$-adaptive FEM strategy described in \cref{sec:numhfem}. We show a sample FEM solution in \cref{fig:exHs}. For the approximation with respect to $z$, we employ $L^2(\Gamma_1)$-SRI and MRI, using $2S-1$ and $S$ uniformly spaced samples of $z\in[10,200]$, respectively, with $S=15$ and $S=25$. This allows us to build rational approximations of the same type ($[S-1]$) with both methods.

The results are shown in \cref{fig:exHf}. There, we can observe that, for a given rational type, $L^2(\Gamma_1)$-SRI and MRI seem to achieve similar approximation errors. Notably, both approaches appear to struggle for low frequencies $z\approx 10$, especially in the cases where fewer snapshots are used.

A comparison of the corresponding runtimes be found in \cref{tab:exHtime}. We can see that, for a given rational type, building the MRI method is more costly than the $L^2(\Gamma_1)$-SRI one, due to the computation of the snapshot Gramian. This is the case despite the lower amount of snapshots required by MRI. Note, in particular, that building the $L^2(\Gamma_1)$-Gramian is considerably faster than computing the $H_0^1(\Omega)$ one, due to the reduced dimensionality of the domain of integration.

\begin{figure}[h!]
\centering
\begin{small}
\begin{tabular}{c|c|c|c|c|}
Method & \multicolumn{2}{c|}{$L^2(\Gamma_1)$-SRI} & \multicolumn{2}{c|}{MRI}\\
\hline
Value of $S$ & 29 & 49 & 15 & 25\\
\hline
\hline
Snapshots & \texttt{2.74e+2} & \texttt{7.74e+2} & \texttt{1.83e+2} & \texttt{3.04e+2}\\
\hline
Build snapshot Gramian & \texttt{1.27e+0} & \texttt{3.21e+0} & \texttt{2.43e+2} & \texttt{7.24e+2}\\
\hline
Assemble surrogate & \texttt{1.42e-2} & \texttt{2.25e-2} & \texttt{9.45e-3} & \texttt{1.18e-2}\\
\hline
\hline
Total & \texttt{2.76e+2} & \texttt{7.77e+2} & \texttt{4.27e+2} & \texttt{1.03e+3}\\
\hline
\end{tabular}
\end{small}
\captionof{table}{Timings in seconds for the construction of $\widetilde{y}$.}
\label{tab:exHtime}
\end{figure}

\subsection{Acoustic scattering of a ``trapping'' cavity}
\label{sec:scattering}

\begin{figure}[p]
\centering
    \begin{minipage}[b]{.325\textwidth}
         \hspace{0.cm}
         \centering
     \begin{subfigure}[t]{.975\textwidth}
         \centering
         \begin{tikzpicture}[scale=3.25]
         \draw (0,0) rectangle (1,1);
         \draw (0.379408389392194,0.343932372542188) -- (0.370591610607806,0.356067627457812) -- (0.613979026519516,0.532898935898998) -- (0.613979026519516,0.632722284235161) -- (0.379408389392194,0.462296740540848) -- (0.370591610607806,0.474431995456472) -- (0.629408389392194,0.662473392204685) -- (0.629408389392194,0.525568004543528) -- (0.379408389392194,0.343932372542188);
         \draw (.125,.15) node {$\partial\Omega'$};
         \draw (.9,.9) node {$\Omega$};
         \draw (.45,.65) node {$\partial D$};
         \end{tikzpicture}
         \caption{$\Omega$}
         \label{fig:exSdomain}
    \end{subfigure}

     \begin{subfigure}[t]{.975\textwidth}
         \centering
         \begin{tikzpicture}[scale=6.5]
         \draw[gray,dashed] (.25,.25) rectangle (.75,.75);
         \draw (0.379408389392194,0.343932372542188) -- (0.370591610607806,0.356067627457812) -- (0.613979026519516,0.532898935898998) -- (0.613979026519516,0.632722284235161) -- (0.379408389392194,0.462296740540848) -- (0.370591610607806,0.474431995456472) -- (0.629408389392194,0.662473392204685) -- (0.629408389392194,0.525568004543528) -- (0.379408389392194,0.343932372542188);
         \draw[very thick, red] (0.370591610607806,0.356067627457812) -- (0.613979026519516,0.532898935898998) -- (0.613979026519516,0.632722284235161) -- (0.379408389392194,0.462296740540848);
         \draw[red] (.5675,.55) node {\large $\omega$};
         \end{tikzpicture}
         \caption{Zoom on the scatterer.}
         \label{fig:exSscatterer}
    \end{subfigure}
    \end{minipage}%
    \begin{minipage}[b]{.3\textwidth}
         \hspace{0.cm}
         \centering
     \begin{subfigure}[t]{.975\textwidth}
         \centering
         \includegraphics[width = .9\textwidth,trim=1.5cm .65cm 1cm .7cm, clip]{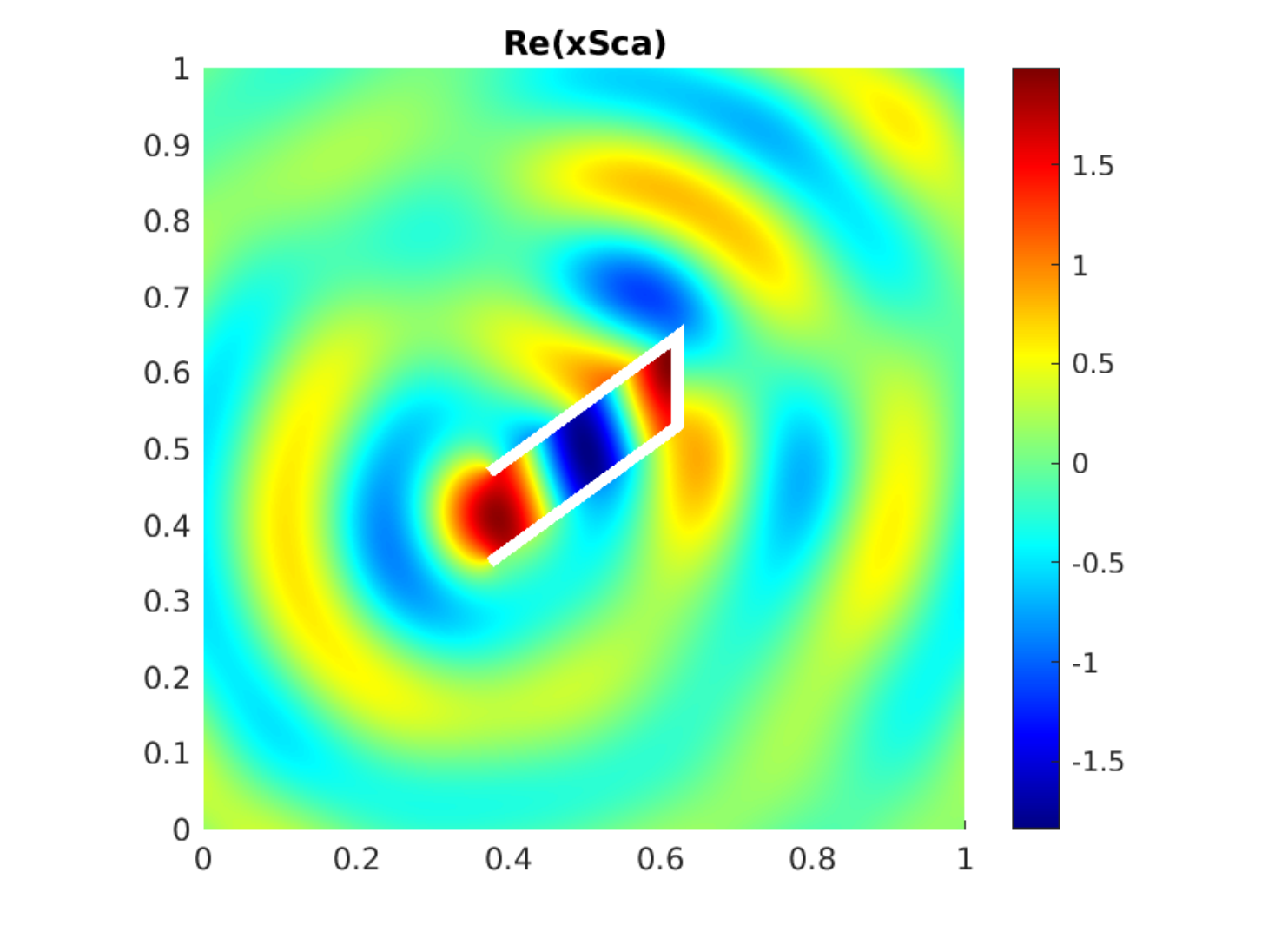}
         \vspace{-.5mm}
         \caption{$\real(u_{89}^\textup{FEM}(25))$}
         \label{fig:exSfemscattered}
    \end{subfigure}

    \begin{subfigure}[t]{.975\textwidth}
         \centering
         \includegraphics[width = .9\textwidth,trim=1.5cm .65cm 1cm .7cm, clip]{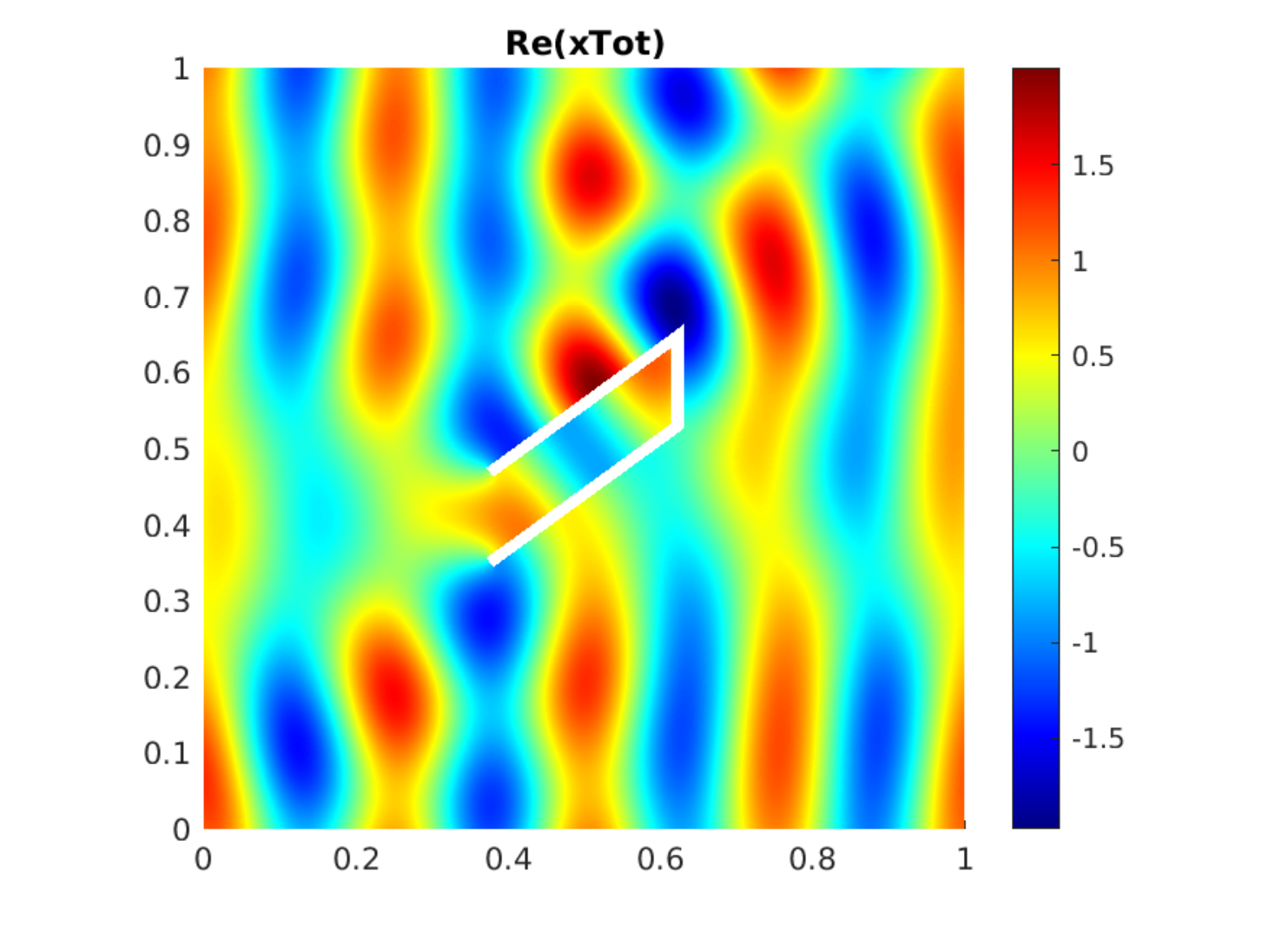}
         \vspace{-.5mm}
         \caption{$\real(u_{89}^\textup{FEM}(25)+u_{\text{inc}}(25))$}
         \label{fig:exSfemtotal}
    \end{subfigure}
    \end{minipage}%
    \begin{minipage}[b]{.3\textwidth}
         \hspace{0.cm}
         \centering
     \begin{subfigure}[t]{.975\textwidth}
         \pgfplotstableread[col sep=comma]{Data/ex7estimator.csv}\tikzdataires
         \centering
         \hspace{-8mm}%
         \begin{tikzpicture}
         \begin{loglogaxis}[
            ticklabel style = {font=\tiny},
            width = \textwidth, height = .8\textwidth,
            xlabel={\small DoFs $N_\bullet$}, ylabel={\small $\eta_\bullet(\Omega)$},
            xlabel style={yshift=.1cm}, ylabel style={yshift=-.25cm},
            ymin = 2e-1, ymax = 9e1,
            legend pos=south west
            ]
            \addplot[semithick, black, mark=*, mark size=.2] plot table [x=dofs25, y=est25]{\tikzdataires};
            \draw[black, densely dashed] (axis cs:1e4,22)--(axis cs:1.21e6,2);
            \draw (axis cs:2e5,15) node[scale=.85, rotate=-35] {$\mathcal{O}(N_\bullet^{-1/2})$};
            \legend{$z=25$};
         \end{loglogaxis}
         \end{tikzpicture}
         \caption{$\eta_\bullet(\Omega)$ vs.\ FEM DoFs.}
         \label{fig:exSr}
     \end{subfigure}

     \begin{subfigure}[t]{.975\textwidth}
         \centering
         \includegraphics[width = .9\textwidth,trim=1.5cm .65cm 1cm .7cm, clip]{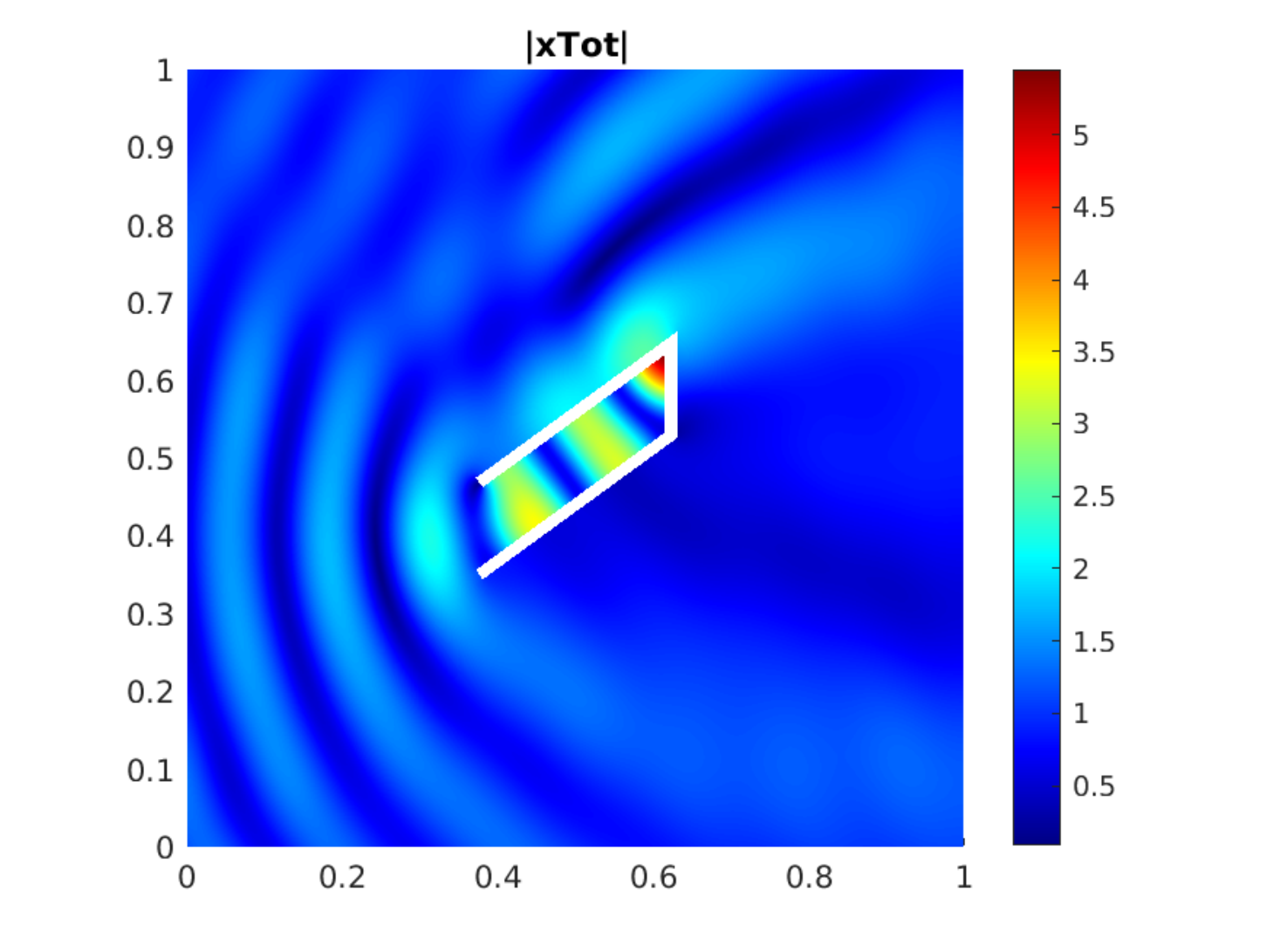}
         \vspace{-.5mm}
         \caption{$|u_{89}^\textup{FEM}(25)+u_{\text{inc}}(25)|$}
         \label{fig:exSfemTotal}
    \end{subfigure}
         \end{minipage}
     \caption{$h$-adaptive FEM results for $z=25$.\vspace{-1mm}}
     \label{fig:exSs}
\end{figure}

We consider a domain of the form $\Omega=\Omega'\setminus D$, where $\Omega'=[0,1]^2$ is a square and $D\subset\Omega'$ is a leftward-facing cavity, in the shape of a slanted ``C'', see \cref{fig:exSdomain}. We are interested in the solution $u=u(z)\in H^1(\Omega)$ of the following Helmholtz equation with impedance boundary conditions (we set $z=k$ here): 
\begin{equation}
\label{eq:exSeqn}
\begin{cases}
    -\Delta u(z)-z^2u(z)=0, &\text{ in } \Omega,\\
    \partial_{\bm{\nu}}u(z)=-\partial_{\bm{\nu}}u_\text{inc}(z), &\text{ on } \partial D,\\
    \partial_{\bm{\nu}}u(z)=\iota zu(z), &\text{ on } \partial\Omega'.
\end{cases}
\end{equation}
The solution $u$ corresponds to the (time-harmonic) wave scattered by the sound-hard scatterer $D$ subject to a horizontal unit plane wave $u_\text{inc}(z)=e^{\iota zx_1}$ incoming from the left. The impedance condition on $\partial\Omega'$ serves as approximation of the Sommerfeld radiation condition at infinity, cf.\ \cref{eq:scattering}. Note that, as in the previous example, the Neumann datum is not affine in $z$. As approximation target, we take the trace $v(z)=u(z)|_\omega$, with $\omega\subset\partial D$ being the trapping boundary of $D$, i.e., the ``interior'' portion of the scatterer's boundary, see \cref{fig:exSscatterer}. Additionally, we consider the (squared) $L^2(\omega)$-norm of $v$, i.e., $y(z)=\int_\omega\abs{v(z)}^2$.

For a given $z$, we compute $u(z)$ via the $h$-adaptive FEM strategy described in \cref{sec:numhfem}, but with $\tol_h=0.5$. We show a sample FEM solution in \cref{fig:exSs}. As in the previous example, we employ $L^2(\omega)$-SRI and MRI for the approximation of $v(z)$ (and consequently $y(z)$) as $z$ varies in the interval of interest $Z=[10,30]$. For both methods, we use $19$ uniformly spaced snapshots, allowing for the construction of rational approximations of types $[9]$ and $[18]$, respectively.

Before proceeding further, we note that, both with $L^2(\omega)$-SRI and with MRI, the approximation of $v(z)\in L^2(\omega)$ will require linear combinations of the snapshot traces $\{v(z_j)\}_{j=1}^{19}$, which live on different discretizations of $\omega$. Still, since $\omega$ is a 1D curve, it is not too expensive to extend all the snapshot traces onto a common discretization of $\omega$, whose vertices are obtained as the union of the locations of all the DoFs of $\{v(z_j)\}_{j=1}^{19}$. See also the discussion in \cref{sec:gsri,sec:mri}.

\begin{figure}[p]
\hspace{0cm}
    \begin{minipage}[t]{.45\textwidth}
    \hspace{0cm}
     \pgfplotstableread[col sep=comma]{Data/ex7global.csv}\tikzdatag
     \begin{tikzpicture}
     \begin{semilogyaxis}[
        scale only axis, grid,
        ticklabel style = {font=\tiny},
        width = .7\textwidth, height = 2.5cm,
        xlabel={\small $z$},
        xlabel style={yshift=.1cm},
        ylabel={$\abs{\widetilde{y}(z)-y(z)}$},
        ylabel style={yshift=-2mm},
        xmin = 10, xmax = 30,
        ymin = 1e-5, ymax = 4,
        name = bottom plot
        ]
        \addplot[blue] plot table [x=k, y=vri_e]{\tikzdatag};
        \addplot[red, densely dashed] plot table [x=k, y=mri_e]{\tikzdatag};
        \addplot[black, only marks, mark = x, mark size = .5mm, domain = 10.001:29.999, samples = 19] {1e0};
     \end{semilogyaxis}

     \begin{axis}[
        scale only axis, grid,
        ticklabel style = {font=\tiny},
        width = .7\textwidth, height = 2.5cm,
        xmin = 10, xmax = 30,
        ylabel={$y(z)$},
        ylabel style={yshift=-2mm},
        xticklabels = {,,},
        ymin = .1, ymax = 9,
    	at={(bottom plot.above north west)},
    	anchor=below south west
        ]
        \addplot[blue] plot table [x=k, y=vri]{\tikzdatag};
        \addplot[red, densely dashed] plot table [x=k, y=mri]{\tikzdatag};
        \addplot[black, only marks, mark = *, mark size = .15mm] plot table [x=k, y=ex]{\tikzdatag};
     \end{axis}
    \end{tikzpicture}
    \caption{Surrogate of $y$ (top) and corresponding error (bottom). In the top plot, high-fidelity $h$-FEM validation points are also included. On top of the bottom plot, we show with crosses the locations of the sample points.\vspace{-1mm}}
    \label{fig:exSf}
    \end{minipage}\hspace{4mm}%
    \begin{minipage}[t]{.45\textwidth}
    \hspace{0cm}
     \pgfplotstableread[col sep=comma]{Data/ex8trace.csv}\tikzdatat
     \centering
     \begin{tikzpicture}
     \begin{axis}[
        scale only axis, grid,
        ticklabel style = {font=\tiny},
        width = .75\textwidth, height = 2.5cm,
        ylabel={$\real(v(25)|_\xi)$},
        ylabel style={yshift=-2mm},
        xticklabels = {,,},
        xmin = 0, xmax = .6906,
        ymin = -2.1, ymax = 3.5,
        legend pos=north west,
        legend columns = 3,
        legend style={nodes={scale=0.6, transform shape}},
        name = top plot
        ]
        \addplot[blue] plot table [x=x, y=trVRIR]{\tikzdatat};
        \addplot[red, densely dashed] plot table [x=x, y=trMRIR]{\tikzdatat};
        \addplot[thick, black, densely dotted] plot table [x=x, y=trExR]{\tikzdatat};
        \draw[black] (axis cs:.30084,-5) -- (axis cs:.30084,5);
        \draw[black] (axis cs:.40067,-5) -- (axis cs:.40067,5);
        \legend{$L^2(\omega)$-SRI,MRI,$h$FEM};
     \end{axis}

     \begin{semilogyaxis}[
        scale only axis, grid,
        ticklabel style = {font=\tiny},
        width = .75\textwidth, height = 2.5cm,
        xlabel={\small $\xi$},
        xlabel style={yshift=.1cm},
        ylabel={\small $\abs{\widetilde{v}(25)|_\xi-v(25)|_\xi}$},
        ylabel style={yshift=-2mm},
        xmin = 0, xmax = .6906,
        ymin = 1e-4, ymax = 1e-1,
    	at={(top plot.below south west)},
    	anchor=north west
        ]
        \addplot[blue] plot table [x=x, y=errVRI]{\tikzdatat};
        \addplot[red, densely dashed] plot table [x=x, y=errMRI]{\tikzdatat};
        \draw[black] (axis cs:.30084,1e-8) -- (axis cs:.30084,1e1);
        \draw[black] (axis cs:.40067,1e-8) -- (axis cs:.40067,1e1);
     \end{semilogyaxis}
    \end{tikzpicture}
    \caption{Surrogate for $v(25)$ (top) and corresponding error (bottom). The high-fidelity $h$-FEM solution is included for validation. We use an (arc-length) parameter $\xi$ to parametrize $\omega$. The two interior vertices of $D$ are denoted by vertical black lines.\vspace{-1mm}}
    \label{fig:exSt}
    \end{minipage}
\end{figure}

\begin{figure}[p]
\centering
\begin{small}
\begin{tabular}{c|c|c|}
Method & $L^2(\omega)$-SRI & MRI\\
\hline
\hline
Snapshots & \multicolumn{2}{c|}{\texttt{3.35e+3}}\\
\hline
Build snapshot Gramian & \texttt{2.75e+0} & \texttt{6.62e+3}\\
\hline
Assemble surrogate & \texttt{1.43e+0} & \texttt{9.49e-1}\\
\hline
\hline
Total & \texttt{3.35e+3} & \texttt{9.97e+3}\\
\hline
\end{tabular}
\end{small}
\captionof{table}{Timings in seconds for the construction of $\widetilde{y}$.}
\label{tab:exStime}
\end{figure}

In \cref{fig:exSf}, we display the approximation of the scalar real-quadratic QoI $y(z)$ over the range of frequencies $z\in Z$. While the $L^2(\omega)$-SRI surrogate approximates the exact $y$ fairly well, we note the appearance of some oscillations in the MRI approximation, especially for $z\approx 30$. This behaviour is due to two compound effects: on one hand, MRI mistakenly places a surrogate pole quite close to the real axis, at $z\approx 29.1-0.14\iota$. On the other hand, modest Runge oscillations can be observed throughout the wavenumber range. We believe that the observed instabilities are due to the ``interpolation'' nature of the $h$-MRI method. Indeed, the employed offline information is a set of snapshots computed via $h$-adaptive FEM, which are affected by the $h$-adaptivity-induced numerical noise $(u_h-u)$. Since the MRI surrogate interpolates the (noisy) snapshots, it is impacted by numerical instabilities. In contrast, the $L^2(\omega)$-SRI relies on a least-squares formulation, which helps in filtering out (part of) the numerical noise due to $h$-adaptivity.


In \cref{fig:exSt}, we show the approximation of $v(z')$ at the arbitrarily chosen frequency $z'=25$, which is \emph{not} one of the wavenumbers that are sampled for the construction of the surrogates. For the two rational surrogates, $\widetilde{v}$ is obtained simply by evaluating the rational approximations at $z'$. Looking at the error, we see that the $L^2(\omega)$-SRI surrogate is about one order of magnitude more accurate than the MRI one.

\section{Conclusions}
\label{sec:conclusions}

In the present paper, we have presented several rational-based spatially adaptive MOR methods for the (non-coercive) parametric-in-frequency Helmholtz PDE endowed with mixed Dirichlet/Neumann/Robin boundary conditions. In particular, we have treated (i) SRI, which is useful in case of scalar-valued QoIs; (ii) $\VV$-SRI, which pertains to $\VV$-valued QoI ($\VV$ being a Hilbert space); (iii) MRI, which provides a surrogate for the solution map $u(z)$ itself.
The offline phase of each of the above-mentioned methods entails the solution of the considered problem for a set of values of the wavenumber. The snapshots are computed by means of spatially adaptive FEM, and, as such, they reside in different discrete spaces, which are adapted to the value of the wavenumber and to the local features of each snapshot.
As a projection-based alternative, we have also considered an $h$-adaptive version of POD.

With the target of comparing the methods, we have performed three numerical experiments.
In the first example, we have observed that building the SRI and MRI surrogates is about $20\%$ faster than POD. The PDEs in the second and third example are both endowed with boundary conditions that depend non-affinely on $z$, preventing a straightforward application of POD. In the second example, the two methods are comparable, in the sense that they achieve similar approximation errors, with the $L^2(\omega)$-SRI method being faster than MRI. In our third example, the $L^2(\omega)$-SRI surrogate displays higher accuracy. Indeed, the MRI surrogate is affected by the presence of spurious poles and by numerical instabilities (Runge oscillations) due to the interplay between the two main features of the method, namely, $h$-adaptivity and interpolation of (noisy) snapshots. In contrast, the $L^2(\omega)$-SRI is more stable, since it produces an approximant in the least-square sense.

We remark that, in this work, the snapshots are computed via an $h$-adaptive algorithm designed to approximate the full state $u$ at optimal rate with respect to the energy norm. Another possibility, particularly appropriate if the functional representing the QoI is fixed, is to consider so-called \emph{goal-oriented} adaptive algorithms, which aim at approximating at optimal rate the QoI directly; see, e.g., \cite{op2001,gs2002}.

Moreover, all the algorithms are applied with sample points that are fixed \emph{a priori}. In contrast, it is often desirable to follow a $z$-\emph{adaptive} approach, where the selection of the sample points is driven by a suitable \emph{a posteriori} estimator (note that this estimator would act on $z$, and not on the spatial domain, as $\eta_\bullet$ does). In the framework of projective MOR methods, this can be achieved with the weak-greedy RB technique. In some cases, the greedy selection of snapshots is also possible in the rational-based setting, e.g., in MRI, by means of the \emph{a posteriori} estimator introduced in~\cite{PN2020,PN2022}. The derivation of novel MOR methodologies that combine $z$- and $h$-adaptivity is subject of a forthcoming work. In particular, with such an approach, we expect to be able to counteract some of the showcased instabilities (mainly, the spurious poles) of the $h$-MRI method.

An interesting further topic to be investigated is the extension of the presented $h$-adaptive MOR methods to the multi-parametric setting, i.e., Helmholtz-like problems where additional parameters (geometry, materials, etc.) are present on top of the wavenumber. This more general framework would be based on the parametric-MOR strategy presented in \cite{PN2020}, and is of practical interest for several important applications, e.g., inverse problems, uncertainty quantification and optimal control problems, where non-$h$-adaptive multi-parametric MOR methods have been mostly employed so far (for instance, we mention~\cite{Baur2011,BP2021,Hess2013,Volkwein2008,Xie2018}).

\section*{Acknowledgments}

F.\ Bonizzoni is member of the INdAM Research group GNCS and acknowledges partial support from the European Research Council ERC under the European Union's Horizon 2020 research and innovation program (Grant agreement No.~865751).
D.\ Pradovera acknowledges partial support from the Swiss National Science Foundation through project 182236.
M.\ Ruggeri acknowledges partial support from the Austrian Science Fund (FWF) through the special research program \emph{Taming complexity in partial differential systems} (grant F65).
The authors would also like to acknowledge the kind hospitality of the Erwin Schr\"odinger International Institute for Mathematics and Physics, where part of this research was carried out.

\appendix

\section{Proof of \cref{lem:tracepspan}}\label{sec:traceinterpolationspan}

In this section, we show that the coefficients of the $\VV$-SRI numerator $P_{[N]}^{\VV\textup{-SRI}}$ belong to the span of the snapshots: $\{p_i\}_{i=0}^N\subset\textup{span}\{v_h(z_j)\}_{j=1}^S$.

\begin{proof}[Proof of \cref{lem:tracepspan}]
Without loss of generality, we assume that $N>0$, and that sample and support points are sorted in such a way that there exists $s\in\{0,\ldots,N+1\}$ such that $z_j=\zeta_{j-1}$ for $j\leq s$, while $\{z_j\}_{j=s+1}^S\cap\{\zeta_i\}_{i=s}^N=\emptyset$.

For all $i=0,\ldots,N$, consider a $\VV$-orthogonal decomposition $p_i=p_i'+p_i''$, with $p_i'\in\textup{span}\{v_h(z_j)\}_{j=1}^S$ and $\langle p_i'',v_h(z_j)\rangle_{\VV}=0$ for all $j=1,\ldots,S$. First, we observe that, for all $i<s$, $p_i=q_iv_h(z_i)$ by \cref{eq:interpolationgammaconstraints}, i.e., $p_i''=0$. So, if $s=N+1$, the claim is automatically proven. Otherwise, it remains to prove that $p_s''=\cdots=p_N''=0$.

By the Pythagorean theorem, \cref{eq:interpolationgammatarget} reads
\begin{equation}\label{eq:apppythagoras}
    \sum_{j=s+1}^S\left\|\sum_{i=0}^N\frac{q_iv_h(z_j)-p_i'}{z_j-\zeta_i}\right\|_{\VV}^2+\sum_{j=s+1}^S\left\|\sum_{i=s}^N\frac{p_i''}{z_j-\zeta_i}\right\|_{\VV}^2.
\end{equation}

For the sake of contradiction, assume that the optimal $p_s'',\ldots,p_N''$ are \emph{not all} equal to 0. By inspection of \cref{eq:apppythagoras}, we see that setting all of them to zero cannot increase the value of the target, i.e., by optimality, the second sum must equal 0. However, this is absurd, since the Cauchy matrix
\begin{equation*}
    \begin{bmatrix}
        (z_{s+1}-\zeta_s)^{-1} & \cdots & (z_{s+1}-\zeta_N)^{-1}\\
        \vdots & \ddots & \vdots\\
        (z_S-\zeta_s)^{-1} & \cdots & (z_S-\zeta_N)^{-1}
    \end{bmatrix}\in\C^{(S-s)\times(N-s+1)}
\end{equation*}
(which has more rows than columns) has full rank \cite{Cauchy}, so that
\begin{equation*}
    \sum_{i=s}^N\frac{p_i''}{z_j-\zeta_i}=0\quad\forall j=s+1,\ldots,S\quad\textup{iff}\quad p_i''=0\quad\forall i=s,\ldots,N.
\end{equation*}
\end{proof}

\section{General algorithm for the $\VV$-SRI method}
\label{sec:traceinterpolation}

In \cref{sec:gsri} we have detailed the $\VV$-SRI algorithm in the particular case where the support points are a subset of the sample points. A similar algorithm may be written even in the general case, where no relation between the set of sample points and support points is assumed.

\begin{lemma}[$\VV$-SRI algorithm]\label{lem:latsri}
    Define the $v$-snapshot Gramian
    \begin{equation}\label{eq:interpolationsnapgram}
        G_h^{(v)}=
        \begin{bmatrix}
        \|v_h(z_1)\|_{\VV}^2 & \cdots & \langle v_h(z_S),v_h(z_1)\rangle_{\VV} \\
        \vdots & \ddots & \vdots\\
        \langle v_h(z_1),v_h(z_S)\rangle_{\VV} & \cdots & \|v_h(z_S)\|_{\VV}^2
        \end{bmatrix}
        \in\C^{S\times S},
    \end{equation}
    whose rank is $T\leq S$. Also, let $s\in\{0,\ldots,N+1\}$ be the cardinality of $\{z_j\}_{j=1}^S\cap\{\zeta_i\}_{i=0}^N$. The corresponding $\VV$-SRI exists and admits the following closed-form representation:
    \begin{itemize}
        \item The coefficients $\mathbf{q}=(q_0,\ldots,q_N)^\top$ of the $\VV$-SRI denominator $Q_{[N]}^{\VV\textup{-SRI}}$ can be found as a (normalized) minimal right singular vector of a matrix $G\in\C^{T(S-s)\times(N+1)}$, i.e.,
        \begin{equation}\label{eq:appclaimq}
            \mathbf{q}=\argmin_{\|\mathbf{q}'\|_{\C^{N+1}}^2=1}\left\|G\mathbf{q}'\ \right\|_{\C^{T(S-s)}}^2.
        \end{equation}
        \item The coefficient matrix $\mathring{P}$, see \cref{eq:tracepexpansion}, may not be unique whenever $T<S$. On the other hand, one possible set of values for it can always be found as
        \begin{equation}\label{eq:appclaimp}
            \mathring{P}_{:i}=H_i\mathbf{q}\quad\textup{for }i=0,\ldots,N.
        \end{equation}
    \end{itemize}
    In the proof, we provide expressions in closed form for $G$ and $H_i$, $i=0,\ldots,N$.
\end{lemma}

\begin{remark}\label{rem:latsrisri}
	\cref{lem:latsri} applies also to SRI, by replacing $v$ with $y$. Note that, in that case, $G_h^{(y)}$ has rank $T=1$ and $G$ has size $(S-s)\times(N+1)$.
\end{remark}

\begin{proof}[Proof of \cref{lem:latsri}]
Without loss of generality, we assume that $N>0$, and that sample and support points are sorted in such a way that $z_j=\zeta_{j-1}$ for $j\leq s$, while $\{z_j\}_{j=s+1}^S\cap\{\zeta_i\}_{i=s}^N=\emptyset$. Also, let $C$ and $\mathring{C}$ be the Cauchy matrices
\begin{equation*}
    C=\begin{bmatrix}
        (z_{s+1}-\zeta_0)^{-1} & \cdots & (z_{s+1}-\zeta_N)^{-1}\\
        \vdots & \ddots & \vdots\\
        (z_S-\zeta_0)^{-1} & \cdots & (z_S-\zeta_N)^{-1}
    \end{bmatrix}\in\C^{(S-s)\times(N+1)}
\end{equation*}
and
\begin{equation*}
    \mathring{C}=\begin{bmatrix}
        (z_{s+1}-\zeta_s)^{-1} & \cdots & (z_{s+1}-\zeta_N)^{-1}\\
        \vdots & \ddots & \vdots\\
        (z_S-\zeta_s)^{-1} & \cdots & (z_S-\zeta_N)^{-1}
    \end{bmatrix}\in\C^{(S-s)\times(N-s+1)},
\end{equation*}
respectively. By our assumptions on $N$ and $S$, $\mathring{C}$ has at least as many rows as columns, so that, in particular, it has full column rank \cite{Cauchy}.

For $i=0,\ldots,s-1$, we have $p_i=q_iv_h(z_{i+1})$ by \cref{eq:interpolationgammaconstraints}, so that, trivially,
\begin{equation*}
    \mathring{P}_{ji}=q_i\delta_{j(i+1)}\quad\textup{for }i=0,\ldots,s-1\textup{ and }j=1,\ldots,S,
\end{equation*}
and part of \cref{eq:appclaimp} follows with
\begin{equation*}
    (H_i)_{ji'}=\delta_{j(i+1)}\delta_{i'i}\quad\textup{for }i=0,\ldots,s-1\textup{, }j=1,\ldots,S\textup{, and }i'=0,\ldots,N.
\end{equation*}
We now move to the cases $i=s,\ldots,N$.

Given $G_h^{(v)}$ in \cref{eq:interpolationsnapgram}, we first compute its rank-revealing Cholesky factorization $G_h^{(v)}=R^HR$, so that
\begin{equation}\label{eq:rightcholesky}
    R=\left[\begin{matrix}
    R_{1:} \\ \vdots\\ R_{T:}
    \end{matrix}\right]=[R_{:1}\cdots R_{:S}]\in\C^{T\times S}.
\end{equation}
Note that, given the matrix $R$, we can find a $\VV$-orthonormal basis $\{\psi_{j'}\}_{j'=1}^T\subset\VV$ satisfying \cref{eq:householder}. For instance, if $T=S$, we can set
\begin{equation*}
\psi_{j'}=\sum_{j=1}^Sv_h(z_j)\left(R^{-1}\right)_{jj'}\qquad\text{for }j'=1,\ldots,S.
\end{equation*}
(Note that such basis is never explicitly needed in the algorithm.)

In particular, by \cref{eq:tracepexpansion,eq:householder}, we have the expansion
\begin{equation}\label{eq:tsrinumorthobasis}
    p_i=\sum_{j=1}^S\sum_{j'=1}^T\mathring{P}_{ji}\psi_{j'}R_{j'j}=\sum_{j'=1}^T(R_{j':}\mathring{P}_{:i})\psi_{j'}\quad\text{for }i=0,\ldots,N.
\end{equation}
This, in turn, allows us to express \cref{eq:interpolationgammatarget} in the Euclidean norm by the Pythagorean theorem
\begin{align}
    \sum_{j=s+1}^S\left\|\sum_{i=0}^N\frac{q_iv_h(z_j)-p_i}{z_j-\zeta_i}\right\|_{\VV}^2=&\sum_{j=s+1}^S\left\|\sum_{j'=1}^T\psi_{j'}\sum_{i=0}^N\frac{q_iR_{j'j}-R_{j':}\mathring{P}_{:i}}{z_j-\zeta_i}\right\|_{\VV}^2\nonumber\\
    =&\sum_{j=s+1}^S\sum_{j'=1}^T\left|\sum_{i=0}^N\frac{q_iR_{j'j}-R_{j':}\mathring{P}_{:i}}{z_j-\zeta_i}\right|^2\label{eq:appinterpolationgammatarget}\\
    =&\sum_{j=s+1}^S\left\|\sum_{i=0}^N\frac{q_iR_{:j}-R\mathring{P}_{:i}}{z_j-\zeta_i}\right\|_{\C^T}^2.\nonumber
\end{align}

It suffices to take the gradient with respect to $\mathring{P}_{:i}$, for $i=s,\ldots,N$, to obtain the optimality conditions characterizing the numerator:
\begin{equation*}
    \mathbf{0}=2\sum_{j=s+1}^S\frac{R^H}{\overline{z_j-\zeta_i}}\sum_{i'=0}^N\frac{R\mathring{P}_{:i'}-q_{i'}R_{:j}}{z_j-\zeta_{i'}}.
\end{equation*}
Using the previously introduced notation, these conditions can be equivalently expressed as
\begin{equation*}
    G_h^{(v)}\sum_{i'=0}^N\left(C^HC\right)_{ii'}\mathring{P}_{:i'}=\sum_{j=1}^{S-s}(C^H)_{ij}(C\mathbf{q})_j(G_h^{(v)})_{:(s+j)}\quad\forall i=s,\ldots,N.
\end{equation*}
Denoting by $\mathbf{e}_{s+j}\in\C^S$ an element of the canonical basis ($(\mathbf{e}_{s+j})_{j'}=\delta_{(s+j)j'}$), we have that $(G_h^{(v)})_{:(s+j)}=G_h^{(v)}\mathbf{e}_{s+j}$, and, for all $i=s,\ldots,N$,
\begin{align}
    G_h^{(v)}\sum_{i'=s}^N\left(C^HC\right)_{ii'}\mathring{P}_{:i'}=&\sum_{j=1}^{S-s}(C^H)_{ij}(C\mathbf{q})_j(G_h^{(v)})_{:(s+j)}-G_h^{(v)}\sum_{i'=0}^{s-1}\left(C^HC\right)_{ii'}\mathring{P}_{:i'}\nonumber\\
    =&G_h^{(v)}\left(\sum_{j=1}^{S-s}(C^H)_{ij}(C\mathbf{q})_j\mathbf{e}_{s+j}-\sum_{i'=0}^{s-1}\left(C^HC\right)_{ii'}H_{i'}\mathbf{q}\right).\label{eq:appinterpolationgammaoptp}
\end{align}

This means that, component-wise,
\begin{equation*}
    \sum_{i'=s}^N\left(C^HC\right)_{ii'}\mathring{P}_{j'i'}=z_{ij'}+\sum_{j=1}^{S-s}(C^H)_{ij}(C\mathbf{q})_j\delta_{(s+j)j'}-\sum_{i'=0}^{s-1}\left(C^HC\right)_{ii'}q_i\delta_{j'(i'+1)}
\end{equation*}
for $i=s,\ldots,N$ and $j'=1,\ldots,S$, with $\mathbf{z}_i=(z_{i1},\ldots,z_{iS})^\top\in\C^S$ some arbitrary element of the kernel of $G_h^{(v)}$.
We observe that the right-hand side above is affine in $\mathbf{q}$, so that we can write
\begin{equation*}
    \sum_{i'=s}^N\left(C^HC\right)_{ii'}\mathring{P}_{j'i'}=z_{ij'}+\mathbf{a}_{ij'}^\top\mathbf{q}\quad\text{for }i=s,\ldots,N,\ j'=1,\ldots,S,
\end{equation*}
with $\mathbf{a}_{ij'}\in\C^{N+1}$, defined entry-wise as
\begin{equation*}
    (\mathbf{a}_{ij'})_{i'}=
    \begin{cases}
        -(C^HC)_{i(j'-1)}&\text{if }j'\leq s\textup{ and }i=i',\\
        0&\text{if }j'\leq s\textup{ and }i\neq i',\\
        (C^H)_{i(j'-s)}(C)_{(j'-s)i'}&\text{if }j'>s.
    \end{cases}
\end{equation*}
Now we collect the equations for a given $j'$ and $i=s,\ldots,N$:
\begin{equation*}
    \mathring{C}^H\mathring{C}\begin{pmatrix}
    \mathring{P}_{j's} \\ \vdots \\ \mathring{P}_{j'N}
    \end{pmatrix}=\begin{pmatrix}
    z_{sj'} \\ \vdots \\ z_{Nj'}
    \end{pmatrix}+\begin{bmatrix}
    \mathbf{a}_{sj'}^\top\\ \vdots \\ \mathbf{a}_{N j'}^\top
    \end{bmatrix}\mathbf{q}\quad\textup{for }j'=1,\ldots,S.
\end{equation*}

By setting
\begin{equation*}
    (\mathring{C}^H\mathring{C})^{-1}=D,\quad\widehat{\mathbf{z}}_{j'}=\begin{pmatrix}
    z_{sj'} \\ \vdots \\ z_{Nj'}
    \end{pmatrix},\quad\textup{and}\quad A_{j'}=\begin{bmatrix}
    \mathbf{a}_{sj'}^\top \\ \vdots \\ \mathbf{a}_{Nj'}^\top
    \end{bmatrix}\quad\textup{for }j'=1,\ldots,S,
\end{equation*}
we can write
\begin{equation*}
    \mathring{P}_{:i}=\begin{bmatrix}
    D_{i:}\widehat{\mathbf{z}}_1 \\ \vdots\\ D_{i:}\widehat{\mathbf{z}}_S
    \end{bmatrix}+
    \begin{bmatrix}
    D_{i:} A_1 \\ \vdots\\ D_{i:} A_S
    \end{bmatrix}\mathbf{q}=:\widetilde{\mathbf{z}}_i+H_i\mathbf{q}\quad\textup{for }i=s,\ldots,N.
\end{equation*}
This, with $\widetilde{\mathbf{z}}_i=\mathbf{0}$, gives \cref{eq:appclaimp} for $i=s,\ldots,N$.

Now, note that, since $\mathbf{z}_i$ belongs to the kernel of $G_h^{(v)}=R^HR$, $\mathbf{z}_i$ belongs to the kernel of $R$ as well. Thus, for all $i=s,\ldots,N$ and $j'=1,\ldots,T$,
\begin{align}
    R_{j':}\mathring{P}_{:i}=&R_{j':}\widetilde{\mathbf{z}}_i+R_{j':} H_i\mathbf{q}
    =\sum_{j=1}^S\sum_{i'=s}^NR_{j'j}D_{ii'} z_{i'j}+R_{j':} H_i\mathbf{q}\nonumber\\
    =&\sum_{i'=s}^ND_{ii'}R_{j':}
    \mathbf{z}_{i'}+R_{j':} H_i\mathbf{q}=\sum_{i'=s}^ND_{ii'}(R
    \mathbf{z}_{i'})_{j'}+R_{j':} H_i\mathbf{q}=R_{j':} H_i\mathbf{q}.\label{eq:appclaimpaffine}
\end{align}
In particular, this means that, for a fixed $\mathbf{q}$, $\mathring{P}_{:i}$ may not be uniquely determined if $G_h^{(v)}$ is rank-deficient (i.e., if $T<S$), even though $p_i$ is always unique, see \cref{eq:tsrinumorthobasis}.

Plugging \cref{eq:appclaimpaffine} into \cref{eq:appinterpolationgammatarget}, after proper re-indexing, yields
\begin{multline*}
    \sum_{j=s+1}^S\sum_{j'=1}^T\left|\sum_{i=0}^N\frac{q_iR_{j'j}-R_{j':} H_i\mathbf{q}}{z_j-\zeta_i}\right|^2\\
    =\sum_{j=s+1}^S\sum_{j'=1}^T\left|\sum_{i=0}^N\frac{q_iR_{j'j}}{z_j-\zeta_i}-\sum_{i,i'=0}^N\sum_{j''=1}^S\frac{R_{j'j''}(H_{i'})_{j''i}q_i}{z_j-\zeta_{i'}}\right|^2\\
    =\sum_{j=s+1}^S\sum_{j'=1}^T\left|\sum_{i=0}^N\left(\frac{R_{j'j}}{z_j-\zeta_i}-\sum_{i'=0}^N\sum_{j''=1}^S\frac{R_{j'j''}(H_{i'})_{j''i}}{z_j-\zeta_{i'}}\right)q_i\right|^2.
\end{multline*}
By setting
\begin{equation*}
g_{ijj'}=\frac{R_{j'j}}{z_j-\zeta_i}-\sum_{i'=0}^N\sum_{j''=1}^S\frac{R_{j'j''} (H_{i'})_{j''i}}{z_j-\zeta_{i'}},    
\end{equation*}
we deduce \cref{eq:appclaimq}, with
\begin{equation*}
    G=\begin{bmatrix}
    g_{0(s+1)1} & \cdots & g_{N(s+1)1}\\
    \vdots & & \vdots\\
    g_{0(s+1)T} & \cdots & g_{N(s+1)T}\\
    g_{0(s+2)1} & \cdots & g_{N(s+2)1}\\
    \vdots & & \vdots\\
    g_{0ST} & \cdots & g_{NST}
    \end{bmatrix}.
\end{equation*}
\end{proof}

\newcommand{\etalchar}[1]{$^{#1}$}

\end{document}